\pdfoutput=1
\RequirePackage{ifpdf}
\ifpdf 
\documentclass[pdftex]{sigma}
\else
\documentclass{sigma}
\fi

\numberwithin{equation}{section}

\newtheorem{Theorem}{Theorem}[section]
\newtheorem*{Theorem*}{Theorem}
\newtheorem{Corollary}[Theorem]{Corollary}

 { \theoremstyle{definition}

\newtheorem{Example}[Theorem]{Example}
\newtheorem{Remark}[Theorem]{Remark} }

\newcommand{\sech}{\operatorname{sech}}

\begin{document}

\renewcommand{\thefootnote}{}

\newcommand{\arXivNumber}{2211.01247}

\renewcommand{\PaperNumber}{015}

\FirstPageHeading

\ShortArticleName{Superposition Formulae for the Geometric B\"acklund Transformations}

\ArticleName{Superposition Formulae for the Geometric B\"acklund\\ Transformations of the Hyperbolic and Elliptic\\ Sine-Gordon and Sinh-Gordon Equations\footnote{This paper is a~contribution to the Special Issue on Symmetry, Invariants, and their Applications in honor of Peter J.~Olver. The~full collection is available at \href{https://www.emis.de/journals/SIGMA/Olver.html}{https://www.emis.de/journals/SIGMA/Olver.html}}}

\Author{Filipe KELMER and Keti TENENBLAT}
\AuthorNameForHeading{F.~Kelmer and K.~Tenenblat}
\Address{Department of Mathematics, Universidade de Bras\'\i lia, Brazil}
\Email{\href{mailto:kelmer.a.f@gmail.com}{kelmer.a.f@gmail.com}, \href{mailto:k.tenenblat@mat.unb.br}{k.tenenblat@mat.unb.br}}

\ArticleDates{Received August 10, 2023, in final form February 03, 2024; Published online February 15, 2024}

\Abstract{We provide superposition formulae for the six cases of B\"acklund transformations corresponding to space-like and time-like surfaces in the 3-dimensional pseudo-Euclidean space. In each case, the surfaces have constant negative or positive Gaussian curvature and they correspond to solutions of one of the following equations: the sine-Gordon, the sinh-Gordon, the elliptic sine-Gordon and the elliptic sinh-Gordon equation. The superposition formulae provide infinitely many solutions algebraically after the first integration of the B\"acklund transformation. Such transformations and the corresponding superposition formulae provide solutions of the same hyperbolic equation, while they show an unusual property for the elliptic equations. The B\"acklund transformation alternates solutions of the elliptic sinh-Gordon equation with those of the elliptic sine-Gordon equation and the superposition formulae provide solutions of the same elliptic equation. Explicit examples and illustrations are given.}


\Classification{58J72; 35J15; 35L10; 53A35; 53A05}

\section{Introduction}
In \cite{KRT2022}, we provided a unified geometric proof
for the six cases of B\"acklund type theorem and integrability
theorem for surfaces in 3-dimensional pseudo-Euclidean space, ${\mathbb{R}}^3_s$, $s=0,1$. These cases include the classical results ($s=0$) and
also those in the Lorentz--Minkowski 3-space ($s=1$), previously obtained by McNertney \cite{McNertney1980}, Palmer \cite{Palmer1990}, Tian \cite{Tian1997} and Gu--Hu--Inoguchi~\cite{Gu--Hu--Inoguchi2002}.

Considering two surfaces $M^2_r,\bar{M}^2_{\bar{r}}\subset{\mathbb{R}}^3_s$ where $r$ and $\bar{r}$ are the indices, $0\leq r,\bar{r}\leq s\leq 1$, the B\"acklund theorem states that a necessary condition for the existence of a B\"acklund-type line congruence in ${\mathbb{R}}^3_s$ between two focal surfaces $M^2_r$ and $\bar{M}^2_{\bar{r}}$, is that both surfaces have the same constant Gaussian curvature, which may be positive or negative depending on~$r$ and~$\bar{r}$. The integrability theorem shows that starting with a space-like or time-like surface ${M\subset \mathbb{R}^3_s}$, there exists a 2-parameter family of surfaces associated to $M$ by a B\"acklund-type line con\-gruence.\looseness=1

Since each such surface corresponds to a solution of one of the following differential equations: sine-Gordon equation, sinh-Gordon equation, elliptic sine-Gordon and elliptic sinh-Gordon, the analytic interpretation of these results provides B\"acklund transformations for the solutions of the corresponding equations. More precisely, such a transformation provides two-parameter new solutions for the sine-Gordon (resp.\ sinh-Gordon) equation from a given solution of the sine-Gordon (resp.\ sinh-Gordon) equation, while the transformation for the elliptic equations provides two-parameter new solutions of the elliptic sinh-Gordon equation from a given solution of the elliptic sine-Gordon equation and vice-versa.

Besides the geometric interpretation of the solutions of these equations, the motivation for studying the (hyperbolic) sine-Gordon, sinh-Gordon and the elliptic sine-Gordon and sinh-Gordon equations comes from their applications in modelling several physical phenomena. The sine-Gordon and sinh-Gordon equations are related to one-dimensional classical field theories, to the theory of crystal dislocations, splay waves in lipid membranes, Bloch--Wall motion, propagation of ultra-short plane wave optical pulses in certain resonant media, etc.\ (see, for example, \cite{AKNS} and
references therein). The elliptic sine-Gordon equation is known to be related to Josephson effects, super-conductors, and spin waves in ferromagnets. The elliptic sinh-Gordon arises as a~model in plasma physics, see, for instance, \cite{FokasLenPelloni, GutLip,JawKaup, Pelloni} and references therein.
This wide range of applications in physics and the fact that there are solutions of the equations that do not correspond to any surface justifies a strictly analytic theory of B\"acklund transformations and superposition formulae for these equations.

In this paper, we provide the analytic superposition formulae for
the B\"acklund transformations discussed in \cite{KRT2022}. The geometric permutability property of the B\"acklund theorem for distinct parameters has an analytic version as a superposition formula, which gives new solutions by algebraic expressions. The importance of a superposition formula is due to the fact that it provides, algebraically, infinitely many new solutions of the partial differential equation, after the first integration required by the B\"acklund transformation.
We observe that the sine-Gordon equation has two B\"acklund transformations \eqref{Backcase1} and \eqref{Backcase2}, which are complementary in the sense that a parameter can be chosen in the intervals $(-1,1)$ or $(1,\infty)$. Similarly, the same occurs with the B\"acklund transformation of the sinh-Gordon equation (see \eqref{Backcase3} and~\eqref{Backcase4}). The superposition formulae hold for distinct parameters chosen in the same interval (see Theorems \ref{SG}--\ref{Sinhgordon}).
 Our superposition formulae include two unusual cases involving the elliptic sine-Gordon and elliptic sinh-Gordon equations
 (see Theorems \ref{elipticSinhgordon} and \ref{elipticSinegordon}). Namely, starting with a solution of the elliptic sine-Gordon (resp.\ elliptic sinh-Gordon) equation, by integrating the B\"acklund transformation for a parameter in the
 interval $(1,\infty)$, we get solutions of the elliptic sinh-Gordon (resp.\ elliptic sine-Gordon) equation and the superposition formula, for distinct parameters, provides algebraically new solutions of the original equation we started with. In the last section of the paper, we provide some examples of solutions for the differential equations and we visualize their graphs and some illustrations of the corresponding surfaces. We also comment on the domain of the solutions.

 An important motivation for a strictly analytic approach to the superposition formulae considered in this paper, is given by the following facts. There is no reason for those interested in obtaining solutions of these equations to go through the geometric setting. Moreover, trying to get the analytic results based on the geometric permutability theorem faces some technical issues, since not all solutions of the differential equations correspond to surfaces. In general, the fundamental forms may present singularities even if the corresponding solution is regular. For instance, the zero solution does not correspond to any surface, since the first fundamental form of such a surface would be singular everywhere. Therefore, starting with such a solution, in order to obtain new solutions of the differential equations, one needs to provide a purely analytic proof of what is called a B\"acklund transformation and a superposition formula, for the solutions of the PDEs. Moreover, even if one is interested in obtaining examples of surfaces, the superposition formulae we provide allow the construction of explicit examples. In particular the new and interesting examples of Section \ref{examples} cannot be obtained solely by applying the geometric theory of surfaces. For instance, in Example~\ref{example2}, one can explicitly obtain surface $X_{12}$ and $X_{123}$ by using the superposition formula. Therefore, knowing the analytic superposition formula for solutions is an important complement to the geometric theory.\looseness=1

We observe that the classical geometric theory of B\"acklund theorem and permutability theorem for surfaces has been extended to
higher-dimensional submanifolds of constant sectional curvature (considering the first and second fundamental forms satisfying Gauss, Codazzi and Ricci equations) in a sequel of papers \cite{TenenblatTerng,Terng} (submanifolds of the Euclidean ambient), \cite{Tenenblat1985} (submanifolds of space forms), \cite{BarbosaFerreiraTenenblat} (Riemannian submanifolds of pseudo-Riemannian space forms), \cite{ChenZuoCheng2004}~(including time-like submanifolds with index $1$ or $n-1$ contained in pseudo-Euclidean space of index $n-1$), see also \cite{DajczerTojeiro}. Considering the first and the second fundamental forms simultaneously diagonalized, the submanifolds are determined by $n\times n$ matrices satisfying a~system of second-order partial differential equations, the so called generalized equations or generating equations in \cite{Tenenblat1998}, The system reduces to the differential equations mentioned above in the 2-dimensional case. Intrinsic geometric generalizations have also been proved, for pseudo-Riemannian manifolds with any index, whose metric has constant sectional curvature~\cite{Campos, CamposTenenblat} see also~\cite{BealsTenenblat1991}. In this case, the metric is given by a vector-valued function satisfying a~system of differential equations the so called intrinsic generalized equations. In \cite[Theorem~2.3]{BarbosaFerreiraTenenblat},
one shows that there is a correspondence between the matrix-valued functions satisfying the generalized equation and the vector-valued function satisfying the intrinsic generalized equation. The generalized equations and the intrinsic generalized equations were renamed as generating equations and intrinsic generating equations in \cite{Tenenblat1998}, since different choices of the index of the metric generate distinct systems of differential equations.

The geometric two-dimensional case is special, in the sense that the first and second fundamental forms of the surfaces are determined by generic functions which satisfy certain partial differential equations. We emphasize that in the 2-dimensional case, the geometric B\"acklund theorem and permutability theorem consider the first and second fundamental forms for surfaces, while the analytic B\"acklund transformation and the superposition formula refer to solutions of the differential equations. Moreover, the geometric results for surfaces have distinct analytic interpretation, according to the differential equation that is being considered. Obtaining and proving superposition formulae for the solutions of the PDEs is not a simple task. In particular, we observe that the geometric permutability Theorem~4.10 between time-like and space-like surfaces obtained in~\cite{ChenZuoCheng2004}
shows the relation between the fundamental forms of the surfaces. Meanwhile, an analytic approach reveals two distinct superposition formulae for the elliptic sine-Gordon equation~\eqref{superposcase6} and for the elliptic sin-Gordon equation \eqref{superposcase5} provided in this paper.

 This paper is organized as follows. In Section~\ref{Prelim}, we recall from \cite{KRT2022}, the six cases of B\"acklund-type theorem (Theorem \ref{TBT}) and the integrability theorem
 (Theorem~\ref{th:integrabilitytheorem}) for surfaces in the 3-dimensional pseudo-Euclidean space ${\mathbb{R}}^3_s$. Moreover, by considering the partial differential equations whose solutions correspond locally to space-like or time-like surfaces with non-zero constant Gaussian curvature (Theorem \ref{PDE}), we provide the analytic interpretation of the geometric results in terms of solutions of these partial differential equations (Theorem \ref{thBT}).
 In Section \ref{Superposition}, we prove the Superposition formulae for the sine-Gordon and the sinh-Gordon equations in Theorem \ref{th:superth1_4}. The superposition formula for the elliptic sinh-Gordon equation is given in Theorem \ref{th:super5} and for the elliptic sine-Gordon equation in Theorem \ref{th:super6}. In Section \ref{sec:6cases},
aiming the applications of the last section and also in order to help the reader to easily access the theory in each case, we give the explicit
B\"acklund transformation and the corresponding superposition formula, for each one of the six cases considered in the previous sections.
 In Section \ref{examples}, we apply the superposition formulae to provide explicit examples and illustrate the theory.

\section{Preliminaries}\label{Prelim}
 We denote by $M_{r}^2\hookrightarrow {\mathbb{R}}_s^3$ a surface of index $r$ isome\-trically immersed in the pseudo-Euclidean space ${\mathbb{R}}^3_{s}$, where the metric in coordinates $(x_1,x_2,x_3)$ is given by ${\rm d}x_1^2+{\rm d}x_2^2+(-1)^s{\rm d}x_3^2$. The indices
$r, s\in \{ 0,1\}$, with $0\leq r\leq s\leq 1$. We consider $X\colon U\subset{\mathbb{R}}^2\rightarrow{\mathbb{R}}_{s}^3$ to be a local parametrization of $M_r^2\subset{\mathbb{R}}_s^3$, where $U$ is an open subset of ${\mathbb{R}}^2$.

 Given a vector $v\in{\mathbb{R}}^3_s$, we say that $v$ is a \textit{unit vector} if the absolute value $|\langle v,v\rangle|=1$.
 We say that $v$ is \textit{space-like} if $\langle v,v\rangle>0$ or $v=0$, \textit{time-like} if $\langle v,v\rangle<0$ and \textit{light-like} if $\langle v,v\rangle=0$ and $v\neq 0$. We say that $M_{r}^2$ is a \textit{space-like} (resp.\ \textit{time-like}) surface if $r=0$ (resp.\ $r=1$).

In \cite{KRT2022}, we introduced the following definition of a B\"acklund-type line congruence in ${\mathbb{R}}^3_s$ for space-like or time-like surfaces in ${\mathbb{R}}^3_s$.
 Let $M^2_{r},\bar{M}^{2}_{\bar{r}}\hookrightarrow{\mathbb{R}}^3_{s}$ be two surfaces that are isometrically immersed in ${\mathbb{R}}_s^3$ with $0\leq r,\bar{r}\leq s\leq 1$.
Let $\psi\colon M_r\rightarrow \bar{M}_{\bar{r}}$ be a diffeomorphism such that for each $p\in M$, $\bar{p}=\psi(p)\neq p$ and
let $v(p)\in {\mathbb{R}}^3_s$ be the vector from $p$ to $\bar{p}$.
We say that $\psi$ is a~{\em B\"acklund-type line congruence} ({\em BLC})
 if there exist constants $\lambda>0$ and $\Lambda\geq 0$, such that
 the following three conditions hold:
 1) The line in ${\mathbb{R}}^3_s$ determined by $v(p)$ is tangent to $M$ and~$\bar{M}$ at~$p$ and $\bar{p}$
respectively.
2) For all $p\in M$, $|v(p)|=\sqrt{|\langle v(p),v(p)\rangle|}=\lambda$.
3) The normals~$N(p)$ and~$\bar{N}(\bar{p})$ are not collinear and $ \bigl\langle N(p),\bar{N}(\bar{p})\bigr\rangle=\Lambda$.
 The constant $\lambda$ is called the {\em distance of the line congruence}.
 We say that the BLC in ${\mathbb{R}}^3_s$ is {\em space-like} (resp.\ {\em time-like}) if for all $p\in M$ the vector $v(p)$ is space-like (resp.\ time-like), i.e., $\epsilon:=\big\langle \frac{v(p)}{\lambda}, \frac{v(p)}{\lambda}\big\rangle=1$ ({resp.} $-1$). One can determine an \textit{angle} between a pair of independent vectors in ${\mathbb{R}}^3_s$ (see \cite{KRT2022}).

 By definition, six cases may occur, according to the values of
 $s$, $r$, $\bar{r}$ and $\epsilon$. We introduce the following {\it notation} that will be used throughout this paper. Consider $\xi\in \{ 1,-1 \}$ and denote
\begin{equation}\label{eq:CS}
C_{\xi}(\phi)=
\begin{cases}
\cos(\phi)&\mbox{if} \ \xi=1,\\
\cosh(\phi)&\mbox{if} \ \xi=-1,
\end{cases}
\qquad \mbox{and}\qquad
S_{\xi}(\phi)=
\begin{cases}
\sin(\phi)&\mbox{if}\ \xi=1,\\
\sinh(\phi)&\mbox{if}\ \xi=-1,
\end{cases}
\end{equation}
where $\phi\in(0,\pi)$ if $\xi=1$ and $\phi\in[0,+\infty)$ if $\xi=-1$. Then $C_{\xi}^2(\phi)+ \xi S_{\xi}^2(\phi)=1$.

 In \cite{KRT2022}, we considered all possible cases of BLCs in ${\mathbb{R}}^3_s$ and provided a unified proof of a~B\"acklund type theorem for surfaces in a pseudo-Euclidean space ${\mathbb{R}}_s^3$, which gives a necessary condition, on the Gaussian curvature for the surfaces. More precisely,

\begin{Theorem}[B\"acklund-type theorem~\cite{KRT2022}]\label{TBT}
Let $M^2_{r},\bar{M}^2_{\bar{r}}\hookrightarrow{\mathbb{R}}^3_{s}$, with $0\leq r,\bar{r}\leq s\leq 1$, be surfaces contained in the pseudo-Euclidean space ${\mathbb{R}}^3_s$. Suppose that $M$ and $\bar{M}$ are related by a B\"acklund type line congruence in ${\mathbb{R}}^3_{s}$. Then both surfaces have the same constant Gaussian curvature, $K=\bar{K}$, given by
\begin{equation}\label{KTB}
K=
\begin{cases}
\displaystyle(-1)^{s+1}\frac{S^2_{\epsilon(-1)^s}(\phi)}{\lambda^2} &\mbox{if} \ r=\bar{r},\vspace{1mm}\\
\displaystyle-\frac{C^2_{-1}(\phi)}{\lambda^2} &\mbox{if} \ r\neq \bar{r},
\end{cases}
\end{equation}
where $\lambda$ is the distance of the line congruence, $\phi$ is the angle between the normals and
 $\epsilon=1$ if the line congruence is space-like or $\epsilon=-1$ if the congruence is time-like. Whenever $\epsilon=-1$, then $s=r=\bar{r}=1$.
\end{Theorem}

In Table~\ref{table1}, we recall from \cite{KRT2022}, the B\"acklund-type theorem for all possible types of BLC.
\begin{table}[th]\renewcommand{\arraystretch}{1.33}
\centering
\begin{tabular}{|l|c|c|c|c|l|}
\hline
\multicolumn{6}{|c|}{\textbf{B\"acklund-type theorem in ${\mathbb{R}}^3_s$}}\\
\hline
 & ${\mathbb{R}}^3_s$ & BLC & 1st surface & 2nd surface & Gaussian curvature\\
\hline
1 & & Euclidian & Euclidian &Euclidian & $K=\bar{K}=-\frac{\sin^2{\phi}}{\lambda^2}$, $\phi\in(0,\pi)$ \\
 &$ s=0$ & $\epsilon=1$ & $r=0$ & $\bar{r}=0$& \\
\hline
2 & &space-like & space-like & space-like & $K=\bar{K}=\frac{\sinh^2{\phi}}{\lambda^2}$, $\phi\in(0,+\infty)$ \\
& $s=1$ &$\epsilon=1$& $r=0$ & $\bar{r}=0$&
\\
\hline
3& & space-like&time-like&time-like&$K=\bar{K}=\frac{\sinh^2{\phi}}{\lambda^2}$, $\phi\in(0,+\infty)$
\\
& $s=1$ &$\epsilon=1$& $r=1$ & $\bar{r}=1$& \\
\hline
4& & time-like&time-like&time-like &$K=\bar{K}=\frac{\sin^2{\phi}}{\lambda^2}$, $\phi\in(0,\pi)$ \\
& $s=1$ &$\epsilon=-1$& $r=1$ & $\bar{r}=1$& \\
\hline
5& &space-like&space-like&time-like&$K=\bar{K}=-\frac{\cosh^2{\phi}}{\lambda^2}$, $\phi\in[0,+\infty)$ \\
& $s=1$ &$\epsilon=1$& $r=0$ & $\bar{r}=1$& \\
\hline
6& &space-like&time-like&space-like&$K=\bar{K}=-\frac{\cosh^2{\phi}}{\lambda^2}$, $\phi\in[0,+\infty)$ \\
& $s=1$ &$\epsilon=1$& $r=1$ & $\bar{r}=0$& \\
\hline
\end{tabular}
\caption{Parameters of the congruence and the curvature for each kind of BLC.}\label{table1}
\end{table}

Moreover, we also proved the integrability theorem, which shows that by starting with a~given surface $M_r$ in ${\mathbb{R}}^3_s$ of nonzero constant Gaussian curvature, one can construct new surfaces~$\bar{M}_{\bar{r}}$ in ${\mathbb{R}}^3_s$, which are locally related to $M$ by a BLC in ${\mathbb{R}}^3_s$.

 From now on, without loss of generality, we can consider $M^2_r\subset{\mathbb{R}}^3_s$ to be a surface with {\it normalized constant Gaussian curvature} $K=\delta\in\{ -1,1\}$. Therefore, the parameters $\phi$ and $\lambda$ are related by the expression (\ref{KTB}), where $K=\delta$.

\begin{Theorem}[integrability theorem \cite{KRT2022}] \label{th:integrabilitytheorem}
Let $M_r^2(\delta)\hookrightarrow{\mathbb{R}}^3_s$ be a surface with constant Gaussian curvature
 $K=\delta=\pm 1$, where $r,s\in\{0,1\}$,
\begin{equation}\label{eq:rs}
0\leq r\leq s= 1 \quad \mbox{if} \quad \delta=1 \qquad \mbox{and}\qquad
0\leq r\leq s\leq 1 \quad \mbox{if}\quad \delta=-1.
\end{equation}
Let $p_0\in M_r^2$ be a non-umbilic point. Given a unit tangent vector ${v}_0\in T_{p_0}M^2_r$ which is not a~principal direction, with $\langle {v}_0,{v}_0\rangle=\epsilon$ and a constant $\phi$, such that
\begin{gather}
\epsilon=\pm 1 \quad\mbox{if} \quad \delta=1 \qquad\mbox{and}\qquad
\epsilon=1 \quad \mbox{if} \quad \delta=-1,
\label{eq:epphi}\\
\phi\in (0,\pi) \quad \mbox{if} \quad \epsilon(-1)^s=1 \qquad\mbox{and}\qquad
\phi\in \left[ 0,+\infty\right) \quad \mbox{if} \quad \epsilon(-1)^s=-1.\nonumber
\end{gather}
Then, there exists a surface $\bar{M}^2_{\bar{r}}\subset{\mathbb{R}}^3_s$, with index $\bar{r}$ satisfying $
(-1)^{\bar{r}}=\delta (-1)^{s+r+1}$,
which is related by a BLC to some open neighborhood of $p_0\in M_r^2$, such that the line of the congruence at~$p_0$ is in the direction of ${v}_0$ and the distance $\lambda$ of the congruence and the inner product $\Lambda$ between the normals at corresponding points
are given by
\begin{equation*}
\lambda=
\begin{cases}
S_{-\epsilon}(\phi)&\mbox{if}\ \delta=1,\\
S_1(\phi)&\mbox{if}\ \delta=-1,\ s=0,\\
C_{-1}(\phi)&\mbox{if}\ \delta=-1, \ s=1,
\end{cases} \qquad
\Lambda=
\begin{cases}
C_{-\epsilon}(\phi)&\mbox{if}\ \delta=1,\\
C_1(\phi)&\mbox{if}\ \delta=-1, \ s=0,\\
S_{-1}(\phi)&\mbox{if}\ \delta=-1, \ s=1,
\end{cases}
\end{equation*}
\noindent where $\phi$ is the angle between the normals.
 In particular, $\lambda$ and $\Lambda$ are related by the relation
\begin{equation}\label{eq:lambdaephi}
\delta \bigl[ (-1)^{s+1}\Lambda^2-\lambda^2\epsilon \bigr]=1.
\end{equation}
\end{Theorem}

\begin{Remark}
{\rm The proof of the integrability theorem shows that given a surface $M_r^2(\delta)\subset{\mathbb{R}}^3_s$, of constant Gaussian curvature $K=\delta=\pm 1$, where $r$, $s$ and $\delta$ satisfy \eqref{eq:rs}, there exists a two-parameter family of surfaces $\bar{M}$ which is related locally to $M$ by BLC in ${\mathbb{R}}^3_s$. The two parameters are $\phi$ (i.e., the constants $\lambda$ and $\Lambda$ related by
\eqref{eq:lambdaephi}) and the one that corresponds to choosing the unit tangent vector ${v}_0$, where
$\langle v_0,v_0\rangle=\epsilon=\pm 1$. Then, Theorem~\ref{TBT} implies that each surface~$\bar{M}$ has Gaussian curvature $\bar{K}=\delta$. Moreover, the proof shows that $\bar{M}$ is obtained by integrating an integrable system.
This procedure is called a \textit{B\"acklund transformation}.
Observe that when $\delta=1$ and $r=1$ then one can choose the unit vector $v_0$ to be spacelike or timelike, i.e., $\epsilon=\pm 1$. Otherwise, if $\delta\neq 1$ or $r\neq 1$, then $v_0$ is always spacelike, i.e., $\epsilon =1$ (see Table~\ref{table1})}.
\end{Remark}

By considering an appropriate local coordinate system for surfaces which admit two real distinct principal curvatures, one can show \cite{KRT2022} that locally a space-like or time-like surface in a pseudo-Euclidean space, with non zero constant Gaussian curvature corresponds to a solution of a partial differential equation.
We observe that there exist time-like surfaces with positive constant Gaussian curvature that are non-umbilic and do not satisfy this condition. For instance, one can have a \textit{time-like totally quasi-umbilic} surface \cite{Clelland2012} (this is, a time-like surface whose shape operator is non-diagonalizable over $\mathbb{C}$) or, a time-like surface whose shape operator is diagonalizable over $\mathbb{C}$ but not over ${\mathbb{R}}$ \cite{Gu--Hu--Inoguchi2002}.

\begin{Theorem}[\cite{KRT2022}] \label{PDE}
Let $M_r^2(\delta)\subset{\mathbb{R}}^3_s$ be a surface with constant Gaussian curvature ${\delta=\pm 1}$, $0\leq r\leq s\leq 1$, such that \eqref{eq:rs} holds. Suppose that $M^2_r$ admits distinct real principal curvatures.
Then, there exists a local parametrization $X(x_1,x_2)$ of $M^2_r$ such that
for $ i,j \in\{1,2\} $,
\[
\left\langle
\frac {X_{x_i}} {\vert X_{x_j}\vert},\frac{X_{x_i}}{\vert X_{x_j}\vert} \right\rangle = (-1)^{(i-1)r} \delta_{ij} \epsilon , \qquad \epsilon=\pm 1,
\]
and a differentiable function $\alpha(x_1,x_2)$ that is a solution of the differential equation
\begin{equation}
\label{eq:PDEgsgordon}
\alpha_{x_1x_1}+\epsilon\delta(-1)^{s} \alpha_{x_2x_2} =-\epsilon\delta S_l(\alpha),\qquad\mbox{where}\quad l=(-1)^{r+s+1}\delta .
\end{equation}
Moreover, in these coordinates, the first and second fundamental forms of $M^2_r$ are given, up to orientation $\tau=\pm 1$, respectively, by
\begin{align}
&{\rm I}=\epsilon C_l^2\Bigl(\frac{\alpha}{2}\Bigr){\rm d}x_1^2+(-1)^r \epsilon S_l^2\Bigl(\frac{\alpha}{2}\Bigr){\rm d}x^2_2,\nonumber\\
&{\rm II}=\tau S_l\Bigl(\frac{\alpha}{2}\Bigr)C_l\Bigl(\frac{\alpha}{2}\Bigr) \bigl[ \epsilon {\rm d}x_1^2-(-1)^r\epsilon l{\rm d}x^2_2\bigr].\label{eq:IandII}
\end{align}
Conversely, for fixed integers $\delta=\pm 1$, $ r, s\in\{0,1\}$, such that \eqref{eq:rs} is satisfied,
 $\epsilon=\pm 1$ and $l=(-1)^{r+s+1}\delta $, let $\alpha(x_1,x_2)$
 be a non zero solution of~\eqref{eq:PDEgsgordon}. Then there exists a surface $M^2_r\subset {\mathbb{R}}^3_s$ with Gaussian curvature $K=\delta$, whose fundamental forms are given by~\eqref{eq:IandII}.
\end{Theorem}

We observe that when $\epsilon\delta(-1)^{s}=-1$, then \eqref{eq:PDEgsgordon} reduces to the sine-Gordon equation if $l=1$, or to the sinh-Gordon equation if $l=-1$, while whenever $\epsilon\delta(-1)^{s}=1$, then \eqref{eq:PDEgsgordon} reduces to the elliptic sine-Gordon equation if $l=1$ or to the elliptic sinh-Gordon equation if $l=-1$,

The analytic version of the integrability theorem, Theorem \ref{th:integrabilitytheorem}, was given in \cite{KRT2022}, by considering the differential equations and their corresponding B\"acklund transformations. 

\begin{Theorem}[\cite{KRT2022}]\label{thBT}
Consider fixed constants $\delta=\pm 1$, $\epsilon=\pm 1$,
$r,s\in \{0,1\}$, $r\leq s$, satisfying~\eqref{eq:rs} and \eqref{eq:epphi}.
Let $\alpha(x_1,x_2)$ be a solution of
the differential equation \eqref{eq:PDEgsgordon}, i.e.,
\begin{equation}\label{eqorig}
\alpha_{x_1x_1}+\epsilon \delta(-1)^{s} \alpha_{x_2x_2} =-\epsilon\delta S_l(\alpha),
\end{equation}
where $l=(-1)^{r+s+1}\delta$ and $S_l$ is given by \eqref{eq:CS}.
Then the following system for $\alpha'(x_1,x_2)$ is integrable:
\begin{gather}
\alpha'_{x_1}-(-1)^s\delta\alpha_{x_2}=\frac{2}{\lambda}S_{(-1)^r}\biggl(\frac{\alpha'}{2}\biggr)C_{l}\Bigl(\frac{\alpha}{2}\Bigr)-\frac{2(-1)^s\tau\Lambda}{\lambda} C_{(-1)^r}\biggl(\frac{\alpha'}{2}\biggr)S_{l}\Bigl(\frac{\alpha}{2}\Bigr),\nonumber\\
\alpha'_{x_2}+\alpha_{x_1}=-\frac{2}{\lambda}C_{(-1)^r}\biggl(\frac{\alpha'}{2}\biggr)S_{l}\Bigl(\frac{\alpha}{2}\Bigr)-\frac{2\delta\tau\Lambda}{\lambda} S_{(-1)^r}\biggl(\frac{\alpha'}{2}\biggr)C_{l}\Bigl(\frac{\alpha}{2}\Bigr),\label{eq:btanalytic_v}
\end{gather}
 where the constants $\lambda$ and $\Lambda$ satisfy \eqref{eq:lambdaephi} and $\tau=\pm 1$. Moreover, $\alpha'$ is a solution of
\begin{equation}\label{eqprime}
\alpha'_{x_1x_1}+\epsilon\delta (-1)^{s} \alpha'_{x_2x_2} =-\epsilon\delta S_{(-1)^r}(\alpha').
\end{equation}
\end{Theorem}

\begin{Remark}\label{rmk:remarkl}
We point out that
 whenever $l=(-1)^r$ \big(i.e., $\delta=(-1)^{s+1}$\big), then \eqref{eqorig} and~\eqref{eqprime} coincide and \eqref{eq:btanalytic_v} is a two parameter (self-)B\"acklund transformation, where one of the parameters is determined by $\lambda$ and $\Lambda$ satisfying \eqref{eq:lambdaephi} and the other one is the initial condition of $\alpha'$ at a point. It provides a procedure of obtaining new solutions of \eqref{eqorig} from a given one, by integrating \eqref{eq:btanalytic_v}.
 When $l\neq (-1)^r $, i.e., $l=(-1)^{r+1}$, then
Theorem \ref{thBT} shows that starting with a solution of \eqref{eqorig} and integrating the two-parameter B\"acklund transformation \eqref{eq:btanalytic_v} it provides solutions of a different equation given by \eqref{eqprime}.
 \end{Remark}

 The analytic B\"acklund transformation \eqref{eq:btanalytic_v} was obtained from a geometric B\"acklund transformation between two surfaces. However, the transformation \eqref{eq:btanalytic_v} holds for every initial solution~$\alpha$ of \eqref{eqorig}, including the trivial one (see Example \ref{example1}), which does not correspond to any surface in ${\mathbb{R}}^3_1$. However, $\alpha'$ does correspond to a surface locally (see Section \ref{examples}).

Although in \cite{KRT2022} we gave a unified proof for Theorem \ref{thBT} in all cases, the superposition formulae for the hyperbolic equations and for the elliptic equations are quite different and require separate proofs (see Theorems \ref{th:superth1_4}--\ref{th:super6}). Therefore,
in view of Remark \ref{rmk:remarkl}, Theorem~\ref{thBT} can be restated in the following two corollaries, considering $l=(-1)^r$, i.e., $\delta=(-1)^{s+1}$, and $l=(-1)^{r+1}$ respectively. The first one will be useful to obtain the superposition formulae for the sine-Gordon and the sinh-Gordon equations, while the second one will be used for the superposition formulae of the elliptic sine-Gordon and elliptic sinh-Gordon equation.

\begin{Corollary}\label{thBT1_4}
Consider fixed constants $\delta=\pm 1$, $\epsilon=\pm 1$,
$r,s\in \{0,1\}$, $r\leq s$, satisfying~\eqref{eq:rs} and~\eqref{eq:epphi} and
$\delta=(-1)^{s+1}$.
Let $\alpha(x_1,x_2)$ be a solution of
the differential equation
\begin{equation}\label{eqorig1_4}
\alpha_{x_1x_1}- \alpha_{x_2x_2} =-\epsilon\delta S_{(-1)^r}(\alpha).
\end{equation}
Then the following system for $\alpha'(x_1,x_2)$ is integrable:
\begin{gather}
\alpha'_{x_1}+\alpha_{x_2}=\frac{2}{\lambda}S_{(-1)^r}\biggl(\frac{\alpha'}{2}\biggr)C_{(-1)^r}\Bigl(\frac{\alpha}{2}\Bigr)-\frac{2(-1)^s\tau\Lambda}{\lambda} C_{(-1)^r}\biggl(\frac{\alpha'}{2}\biggr)S_{(-1)^r}\Bigl(\frac{\alpha}{2}\Bigr),\nonumber\\
\alpha'_{x_2}+\alpha_{x_1}=-\frac{2}{\lambda}C_{(-1)^r}\biggl(\frac{\alpha'}{2}\biggr)S_{(-1)^r}\Bigl(\frac{\alpha}{2}\Bigr)-\frac{2\delta\tau\Lambda}{\lambda} S_{(-1)^r}\biggl(\frac{\alpha'}{2}\biggr)C_{(-1)^r}\Bigl(\frac{\alpha}{2}\Bigr),\label{eq:btanalytic_vc1_4}
\end{gather}
 where $\tau=\pm 1$, $\lambda=S_{-\epsilon\delta}(\phi)$, $\Lambda=C_{-\epsilon\delta}(\phi)$, for a constant $\phi$. Moreover,
 $\alpha'$ is a solution of~\eqref{eqorig1_4}.
\end{Corollary}

For $\alpha$ and $\alpha'$ as in Corollary \ref{thBT1_4}, we will say that $\alpha'$ is associated to $\alpha$ by the {\it $\phi$-B\"acklund transformation \eqref{eq:btanalytic_vc1_4}}.
 The B\"acklund transformation above includes the classical case ($s=0$) and the case obtained in \cite{McNertney1980} when $s=1$, $\epsilon=-1$, $r=\delta=1$ and $\tau=-1$.
 Moreover, it shows that for the sine-Gordon equation (resp.\ the sinh-Gordon equation)
 the transformation \eqref{eq:btanalytic_vc1_4}
 can be used in two complementary ways with respect to the parameters, namely
 $\Lambda \in (-1,1)$ or $\Lambda\in(1,+\infty)$, as one can see in Cases 1 and 2 (resp.\ Cases 3 and 4) of Section \ref{sec:6cases}.

In order to state the second corollary, considering $l=(-1)^{r+1}$ in Theorem \ref{thBT}, we observe that since $l=(-1)^{r+s+1}\delta$, it follows that $\delta=(-1)^s$ and hence \eqref{eq:rs} implies that $s=1$. Therefore, $\delta=-1$ and hence it follows from \eqref{eq:epphi} that $\epsilon=1$.

 \begin{Corollary}\label{thBT5_6}
Consider a fixed constant $r\in \{0,1\}$.
Let $\alpha(x_1,x_2)$ be a solution of
the differential equation
\begin{equation}\label{eqorig5_6}
\alpha_{x_1x_1}+\alpha_{x_2x_2} = S_{(-1)^{r+1}}(\alpha).
\end{equation}
Then the following system for $\alpha'(x_1,x_2)$ is integrable:
\begin{gather}
\alpha'_{x_1}-\alpha_{x_2}=\frac{2}{\lambda}S_{(-1)^r}\biggl(\frac{\alpha'}{2}\biggr)C_{(-1)^{r+1}}\Bigl(\frac{\alpha}{2}\Bigr)+\frac{2\tau\Lambda}{\lambda} C_{(-1)^r}\biggl(\frac{\alpha'}{2}\biggr)S_{(-1)^{r+1}}\Bigl(\frac{\alpha}{2}\Bigr),\nonumber\\
\alpha'_{x_2}+\alpha_{x_1}=-\frac{2}{\lambda}C_{(-1)^r}\biggl(\frac{\alpha'}{2}\biggr)S_{(-1)^{r+1}}\Bigl(\frac{\alpha}{2}\Bigr)+\frac{2\tau\Lambda}{\lambda} S_{(-1)^r}\biggl(\frac{\alpha'}{2}\biggr)C_{(-1)^{r+1}}\Bigl(\frac{\alpha}{2}\Bigr),\label{eq:btanalytic_v5_6}
\end{gather}
 where $\tau=\pm 1$, $\lambda=\cosh\phi$, $\Lambda=\sinh \phi$, $\phi\in [0,\infty)$.
Moreover, $\alpha'$ is a solution of
\begin{equation}\label{eqprime5_6}
\alpha'_{x_1x_1}+\alpha'_{x_2x_2} = S_{(-1)^r}(\alpha').
\end{equation}
\end{Corollary}

For $\alpha$ and $\alpha'$ as in Corollary \ref{thBT5_6}, we will say that $\alpha'$ is associated to $\alpha$ by the {\it $\phi$-B\"acklund transformation \eqref{eq:btanalytic_v5_6}}. We observe that when $r=0$
 this transformation, with $\tau=-1$, was
obtained by Tian \cite{Tian1997}.

\section{Superposition formulae} \label{Superposition}
In this section, we provide the superposition formulae for the partial differential equations considered in the previous section. Such a formula provides algebraically infinitely many new solutions of the partial differential equations, after the first integration required by the corresponding B\"acklund transformation.

Our first result provides the superposition formula for the B\"acklund transformations of the differential equations given in Corollary \ref{thBT1_4}. Then, we prove two theorems giving the superposition formulae for the equations given in Corollary \ref{thBT5_6}.
These two results are unusual and they show that starting with a solution of the elliptic sinh-Gordon equation and applying twice the B\"acklund transformation, with distinct parameters, one gets solutions of the same equation although the intermediate step gives solutions of the elliptic sine-Gordon equation. A similar property occurs starting with a solution of
 the elliptic sine-Gordon equation.

\begin{Theorem}\label{th:superth1_4}
 Consider fixed constants $\delta=\pm 1$, $\epsilon=\pm 1$,
$r,s\in \{0,1\}$, $r\leq s$, satisfying~\eqref{eq:rs} and \eqref{eq:epphi}, and $\delta=(-1)^{s+1}$. Let $\alpha(x_1,x_2)$ be a solution of \eqref{eqorig1_4} and
 let $\alpha'$ $($resp.\ $\alpha'')$ be a~solution of \eqref{eqorig1_4}
 associated to $\alpha$ by the $\phi_1($resp.\ $\phi_2)$-B\"acklund transformation \eqref{eq:btanalytic_vc1_4}.
 If $\phi_1\neq \phi_2$, then there exists a unique solution $\alpha^*$ of
 \eqref{eqorig1_4} such that $\alpha^*$ is associated to $\alpha'$ $($resp.~$\alpha'')$ by the $\phi_2($resp.\ $\phi_1)$-B\"acklund transformation \eqref{eq:btanalytic_vc1_4}. Moreover, $\alpha^*$ is
given algebraically by
\begin{equation}\label{eq:superpos}
T_{(-1)^r}\biggl(\frac{\alpha^*-\alpha}{4}\biggr) =
\delta\tau \frac{S_{-\delta\epsilon}\bigl(\frac{\phi_2+\phi_1}{2}\bigr) } {S_{-\delta\epsilon}\bigl(\frac{\phi_2-\phi_1}{2}\bigr)}
T_{(-1)^r}\biggl(\frac{\alpha'-\alpha''}{4}\biggr),
\end{equation}
where $T_{(-1)^r}$ is the function $\tan$ $($resp.\ $\tanh)$ when $r=0$ $($resp.\ $r=1)$.
\end{Theorem}

\begin{proof}
Let $\alpha(x_1,x_2)$ be a solution of \eqref{eqorig1_4}. Let $\alpha'$ (resp.\ $\alpha''$) be a solution of \eqref{eqorig1_4} associated to $\alpha$ by the B\"acklund transformation \eqref{eq:btanalytic_vc1_4}
for parameters $\{\lambda_1=S_{-\epsilon\delta}(\phi_1),\; \Lambda_1=C_{-\epsilon\delta}(\phi_1)\}$ (resp.\ $\{\lambda_2=S_{-\epsilon\delta}(\phi_2),\; \Lambda_2=C_{-\epsilon\delta}(\phi_2)\}$), i.e., $\alpha'$ and $\alpha''$ satisfy
\begin{gather}
\alpha'_{x_1}+\alpha_{x_2}=\frac{2}{\lambda_1}S_{l}\biggl(\frac{\alpha'}{2}\biggr)C_{l}\Bigl(\frac{\alpha}{2}\Bigr)-\frac{2(-1)^s\tau\Lambda_1}{\lambda_1} C_{l}\biggl(\frac{\alpha'}{2}\biggr)S_{l}\Bigl(\frac{\alpha}{2}\Bigr),\nonumber\\[1mm]
\alpha'_{x_2}+\alpha_{x_1}=-\frac{2}{\lambda_1}C_{l}\biggl(\frac{\alpha'}{2}\biggr)S_{l}\Bigl(\frac{\alpha}{2}\Bigr)-\frac{2\delta\tau\Lambda_1}{\lambda_1} S_{l}\biggl(\frac{\alpha'}{2}\biggr)C_{l}\Bigl(\frac{\alpha}{2}\Bigr),\label{eq:btanalytic_v1}
\\[1mm]
\alpha''_{x_1}+\alpha_{x_2}=\frac{2}{\lambda_2}S_{l}\biggl(\frac{\alpha''}{2}\biggr)C_{l}\Bigl(\frac{\alpha}{2}\Bigr) -\frac{2(-1)^s\tau\Lambda_2}{\lambda_2}C_{l}\biggl(\frac{\alpha''}{2}\biggr)S_{l}\Bigl(\frac{\alpha}{2}\Bigr),\nonumber\\[1mm]
\alpha''_{x_2}+\alpha_{x_1}=-\frac{2}{\lambda_2}C_{l}\biggl(\frac{\alpha''}
{2 }\biggr)S_{l}\Bigl(\frac{\alpha}{2}\Bigr)-\frac{2\delta\tau\Lambda_2}{\lambda_2} S_{l}\biggl(\frac{\alpha''}{2}\biggr)C_{l}\Bigl(\frac{\alpha}{2}\Bigr),\label{eq:btanalytic_v2}
\end{gather}
 where $ l=(-1)^r$.

Assume there exists $\alpha^*$ associated to $\alpha'$ (resp.\ $\alpha''$) by the B\"acklund transformation with parameters
 $\{\lambda_2,\; \Lambda_2\}$ (resp.\ $\{\lambda_1,\; \Lambda_1\}$). then $\alpha^*$ should satisfy the following systems of equations
\begin{gather}
\alpha^{*}_{x_1}+\alpha'_{x_2} = \frac{2}{\lambda_2}S_{l}\biggl(\frac{\alpha^{*}}{2}\biggr)C_{l}\biggl(\frac{\alpha'}{2}\biggr)-\frac{2(-1)^s\tau\Lambda_2}{\lambda_2} C_{l}\biggl(\frac{\alpha^{*}}{2}\biggr)S_{l}\biggl(\frac{\alpha'}{2}\biggr),\nonumber\\
\alpha^{*}_{x_2}+\alpha'_{x_1}=-\frac{2}{\lambda_2}C_{l}\biggl(\frac{\alpha^{*}}{2}\biggr)S_{l}\biggl(\frac{\alpha'}{2}\biggr)-\frac{2\delta\tau\Lambda_2}{\lambda_2} S_{l}\biggl(\frac{\alpha^{*}}{2}\biggr)C_{l}\biggl(\frac{\alpha'}{2}\biggr),\label{eq:btanalytic_v2*}\\
\alpha^*_{x_1}+\alpha''_{x_2}=\frac{2}{\lambda_1}S_{l}\biggl(\frac{\alpha^*}{2}\biggr)C_{l}\biggl(\frac{\alpha''}{2}\biggr)-\frac{2(-1)^s\tau\Lambda_1}{\lambda_1} C_{l}\biggl(\frac{\alpha^*}{2}\biggr)S_{l}\biggl(\frac{\alpha''}{2}\biggr),\nonumber\\
\alpha^*_{x_2}+\alpha''_{x_1}=-\frac{2}{\lambda_1}C_{l}\biggl(\frac{\alpha^*}{2}\biggr)S_{l}\biggl(\frac{\alpha''}{2}\biggr)-\frac{2\delta\tau\Lambda_1}{\lambda_1} S_{l}\biggl(\frac{\alpha^*}{2}\biggr)C_{l}\biggl(\frac{\alpha''}{2}\biggr).\label{eq:btanalytic_v1*}
\end{gather}

Subtracting the first (resp.\ second) equation of \eqref{eq:btanalytic_v2} from the first
(resp.\ second) equation of~\eqref{eq:btanalytic_v1}, we obtain
\begin{gather}
\alpha'_{x_1}-\alpha''_{x_1}= \frac{2}{\lambda_1}S_{l}\biggl(\frac{\alpha'}{2}\biggr)C_{l}\Bigl(\frac{\alpha}{2}\Bigr)-\frac{2(-1)^s\tau\Lambda_1}{\lambda_1} C_{l}\biggl(\frac{\alpha'}{2}\biggr)S_{l}\Bigl(\frac{\alpha}{2}\Bigr) -\frac{2}{\lambda_2}S_{l}\biggl(\frac{\alpha''}{2}\biggr)C_{l}\Bigl(\frac{\alpha}{2}\Bigr)
\nonumber\\ \hphantom{\alpha'_{x_1}-\alpha''_{x_1}=}{}
+\frac{2(-1)^s\tau\Lambda_2}{\lambda_2}C_{l}
\biggl(\frac{\alpha''}{2}\biggr)S_{l}\Bigl(\frac{\alpha}{2}\Bigr),\label{eq:primv1-v2}
\\
\alpha'_{x_2}-\alpha''_{x_2}=
-\frac{2}{\lambda_1}C_{l}\biggl(\frac{\alpha'}{2}\biggr)S_{l}\Bigl(\frac{\alpha}{2}\Bigr)-\frac{2\delta\tau\Lambda_1}{\lambda_1} S_{l}\biggl(\frac{\alpha'}{2}\biggr)C_{l}\Bigl(\frac{\alpha}{2}\Bigr) +\frac{2}{\lambda_2}C_{l}\biggl(\frac{\alpha''}
{2 }\biggr)S_{l}\Bigl(\frac{\alpha}{2}\Bigr)
\nonumber\\ \hphantom{\alpha'_{x_2}-\alpha''_{x_2}=}{}
+\frac{2\delta\tau\Lambda_2}{\lambda_2} S_{l}\biggl(\frac{\alpha''}{2}\biggr)C_{l}\Bigl(\frac{\alpha}{2}\Bigr).\label{eq:segv1-v2}
\end{gather}

Similarly, subtracting the first (resp.\ second) equation of \eqref{eq:btanalytic_v1*} from the first (resp.\ second) equation of \eqref{eq:btanalytic_v2*}, we obtain
\begin{gather}
\alpha'_{x_2}-\alpha''_{x_2} = \frac{2}{\lambda_2}S_{l}\biggl(\frac{\alpha^*}{2}\biggr)C_{l}\biggl(\frac{\alpha'}{2}\biggr)+ \frac{2\delta\tau\Lambda_2}{\lambda_2} C_{l}\biggl(\frac{\alpha^*}{2}\biggr)S_{l}\biggl(\frac{\alpha'}{2}\biggr)-\frac{2}{\lambda_1}S_{l}\biggl(\frac{\alpha^*}{2}\biggr)C_{l}\biggl(\frac{\alpha''}{2}\biggr)
\nonumber\\ \hphantom{\alpha'_{x_2}-\alpha''_{x_2} =}{}
-\frac{2\delta\tau\Lambda_1}{\lambda_1} C_{l}\biggl(\frac{\alpha^*}{2}\biggr)S_{l}\biggl(\frac{\alpha''}{2}\biggr),\label{eq:primv2*-v1*}
\\
\alpha'_{x_1}-\alpha''_{x_1}=
-\frac{2}{\lambda_2}C_{l}\biggl(\frac{\alpha^*}{2}\biggr)S_{l}\biggl(\frac{\alpha'}{2}\biggr)-\frac{2\delta\tau\Lambda_2}{\lambda_2} S_{l}\biggl(\frac{\alpha^*}{2}\biggr)C_{l}\biggl(\frac{\alpha'}{2}\biggr)+\frac{2}{\lambda_1}C_{l}\biggl(\frac{\alpha^*}{2}\biggr)S_{l}\biggl(\frac{\alpha''}{2}\biggr)
\nonumber\\ \hphantom{\alpha'_{x_1}-\alpha''_{x_1}=}{}
+\frac{2\delta\tau\Lambda_1}{\lambda_1} S_{l}\biggl(\frac{\alpha^*}{2}\biggr)C_{l}\biggl(\frac{\alpha''}{2}\biggr). \label{eq:segv2*-v1*}
\end{gather}
Now comparing equations \eqref{eq:primv1-v2} with \eqref{eq:segv2*-v1*} and multiplying all terms by $\lambda_1\lambda_2$, we get
\begin{gather*} 
\lambda_1
\biggl[ C_l\biggl(\frac{\alpha^*}{2}\biggr) S_l\biggl(\frac{\alpha'}{2}\biggr)
-S_l\biggl(\frac{\alpha''}{2}\biggr) C_l\Bigl(\frac{\alpha}{2}\Bigr)\biggr]\\
\qquad{}+\lambda_1\Lambda_2\delta\tau
\biggl[ S_l\biggl(\frac{\alpha^*}{2}\biggr) C_l\biggl(\frac{\alpha'}{2}\biggr)
-C_l\biggl(\frac{\alpha''}{2}\biggr) S_l\Bigl(\frac{\alpha}{2}\Bigr)\biggr] \\
\qquad{}+\lambda_2 \biggl[ S_l\biggl(\frac{\alpha'}{2}\biggr) C_l\Bigl(\frac{\alpha}{2}\Bigl)
 - C_l\biggl(\frac{\alpha^*}{2}\biggr) S_l\biggl(\frac{\alpha''}{2}\biggr) \biggr]\\
 \qquad{}+\lambda_2\Lambda_1\delta\tau \biggl[ C_l\biggl(\frac{\alpha'}{2}\biggr) S_l\Bigl(\frac{\alpha}{2}\Bigr)
 - S_l\biggl(\frac{\alpha^*}{2}\biggr) C_l\biggl(\frac{\alpha''}{2}\biggr) \biggr]=0,
\end{gather*}
where we have used the fact that $\delta=(-1)^{s+1}$.

Using trigonometric identities,
the above equation reduces to
\begin{gather*}
\biggl[S_l\biggl(\frac{\alpha^*+\alpha'}{2}\biggr) - S_l\biggl(\frac{\alpha''+\alpha}{2}\biggr) \biggr]
(\lambda_1 +\delta\tau\lambda_1\Lambda_2)\\
 \qquad{} +\biggl[S_l\biggl(\frac{\alpha^*-\alpha'}{2}\biggr) - S_l\biggl(\frac{\alpha+\alpha''}{2}\biggr) \biggr]
(-\lambda_1 +\delta\tau\lambda_1\Lambda_2) \\
 \qquad{} +\biggl[S_l\biggl(\frac{\alpha+\alpha'}{2}\biggr) - S_l\biggl(\frac{\alpha^*+\alpha'}{2}\biggr) \biggr]
(\lambda_2 +\delta\tau\lambda_2\Lambda_1)\\
 \qquad{} +\biggl[S_l\biggl(\frac{\alpha-\alpha'}{2}\biggr) - S_l\biggl(\frac{\alpha^*-\alpha''}{2}\biggr) \biggr]
(-\lambda_2 +\delta\tau\lambda_2\Lambda_1)=0,
\end{gather*}
 which can be rewritten as
 \begin{gather*}
 C_l\biggl(\frac{\alpha^*+\alpha+\alpha'+\alpha''}{4}\biggr) S_l\biggl(\frac{\alpha^*-\alpha+\alpha'-\alpha''}{4}\biggr) (\lambda_1 +\delta\tau\lambda_1\Lambda_2)\\
 \qquad{} + C_l\biggl(\frac{\alpha^*+\alpha-(\alpha'+\alpha'')}{4}\biggr) S_l\biggl(\frac{\alpha^*-\alpha-(\alpha'-\alpha'')}{4}\biggr) (-\lambda_1 +\delta\tau\lambda_1\Lambda_2)\\
 \qquad{} + C_l\biggl(\frac{\alpha^*+\alpha+\alpha'+\alpha''}{4}\biggr) S_l\biggl(\frac{-(\alpha^*-\alpha)+\alpha'-\alpha''}{4}\biggr) (\lambda_2 +\delta\tau\lambda_2\Lambda_1)\\
 \qquad{} + C_l\biggl(\frac{\alpha^*+\alpha-(\alpha'+\alpha'')}{4}\biggr) S_l\biggl(\frac{-(\alpha^*-\alpha)-(\alpha'-\alpha'')}{4}\biggr) (-\lambda_2 +\delta\tau\lambda_2\Lambda_1)=0.
 \end{gather*}

 By applying again trigonometric identities
 to each term $C_l$ and $S_l$ above, we obtain the equation in terms of $S_l$ and $C_l$ of the expressions
\[
\frac{\alpha^*-\alpha}{4}, \qquad
\frac{\alpha'-\alpha''}{4},
\qquad
\frac{\alpha^*+\alpha}{4}, \qquad
\frac{\alpha'+\alpha''}{4}.
\]
After cancellations it reduces to
\begin{gather*}
\delta\tau C_l\biggl(\frac{\alpha^*+\alpha}{4}\biggr) C_l\biggl(\frac{\alpha'+\alpha''}{4}\biggr)
S_l\biggl(\frac{\alpha^*-\alpha}{4}\biggr)C_l\biggl(\frac{\alpha'-\alpha''}{4}\biggr)(\lambda_1\Lambda_2-\lambda_2\Lambda_1)\\
\qquad{} -l \delta\tau S_l\biggl(\frac{\alpha^*+\alpha}{4}\biggr) S_l\biggl(\frac{\alpha'+\alpha''}{4}\biggr)
C_l\biggl(\frac{\alpha^*-\alpha}{4}\biggr)S_l\biggl(\frac{\alpha'-\alpha''}{4}\biggr)(\lambda_1\Lambda_2+\lambda_2\Lambda_1)\\
\qquad{} + C_l\biggl(\frac{\alpha^*+\alpha}{4}\biggr) C_l\biggl(\frac{\alpha'+\alpha''}{4}\biggr)
C_l\biggl(\frac{\alpha^*-\alpha}{4}\biggr)S_l\biggl(\frac{\alpha'-\alpha''}{4}\biggr) (\lambda_1+\lambda_2)\\
\qquad{} -lS_l\biggl(\frac{\alpha^*+\alpha}{4}\biggr) S_l\biggl(\frac{\alpha'+\alpha''}{4}\biggr)
S_l\biggl(\frac{\alpha^*-\alpha}{4}\biggr)C_l\biggl(\frac{\alpha'-\alpha''}{4}\biggr) (\lambda_1-\lambda_2)=0.
\end{gather*}

Dividing by
\[
C_l\biggl(\frac{\alpha^*+\alpha}{4}\biggr) C_l\biggl(\frac{\alpha'+\alpha''}{4}\biggr)C_l\biggl(\frac{\alpha^*-\alpha}{4}\biggr)C_l\biggl(\frac{\alpha'-\alpha''}{4}\biggr),
\]
and introducing the notation $T_l=\displaystyle \frac{S_l}{C_l}$, the equation above reduces to
\begin{gather}
 \delta\tau T_l\biggl(\frac{\alpha^*-\alpha}{4}\biggr)(\lambda_1\Lambda_2-\lambda_2\Lambda_1)
 -l\delta\tau T_l\biggl(\frac{\alpha^*+\alpha}{4}\biggr)T_l\biggl( \frac{\alpha'+\alpha''}{4}\biggr)
T_l\biggl(\frac{\alpha'-\alpha''}{4}\biggr) (\lambda_1\Lambda_2+\lambda_2\Lambda_1) \nonumber\\
\qquad{} +T_l\biggl(\frac{\alpha'-\alpha''}{4}\biggr)(\lambda_1+\lambda_2) \nonumber\\
\qquad{} - lT_l\biggl(\frac{\alpha^*+\alpha}{4} \biggr)T_l\biggl(\frac{\alpha'+\alpha''}{4}\biggr)T_l\biggl(\frac{\alpha^*-\alpha}{4}\biggr)(\lambda_1-\lambda_2)=0.\label{Tl}
\end{gather}

Since
\[
\Lambda_i=C_{-\delta\epsilon}(\phi_i) , \qquad \lambda_i=S_{-\delta\epsilon}(\phi_i), \qquad i=1,2,
\]
we have that
\begin{gather*}
\lambda_1\Lambda_2\pm\lambda_2\Lambda_1=2
S_{-\delta\epsilon}\biggl(\frac{\phi_1\pm \phi_2}{2}\biggr)C_{-\delta\epsilon}\biggl(\frac{\phi_1\pm \phi_2}{2}\biggr),
\\
\lambda_1\pm\lambda_2= 2 S_{-\delta\epsilon}\biggl(\frac{\phi_1\pm \phi_2}{2}\biggr)
C_{-\delta\epsilon}\biggl(\frac{\phi_1\mp \phi_2}{2}\biggr).
\end{gather*}
Substituting these expressions into \eqref{Tl}, we conclude that the following product vanishes
\begin{gather*}
 \biggl[T_l\biggl(\frac{\alpha^*-\alpha}{4}\biggr) S_{-\delta\epsilon}\biggl(\frac{\phi_2-\phi_1}{2}\biggr)-\delta\tau
T_l\biggl(\frac{\alpha'-\alpha''}{4}\biggr) S_{-\delta\epsilon}\biggl(\frac{\phi_2+\phi_1}{2}\biggr)\biggr] \\
 \qquad{}\times\biggl[ -\delta\tau C_{-\delta\epsilon}\biggl(\frac{\phi_1-\phi_2}{2}\biggr)+l T_l\biggl(\frac{\alpha^*+\alpha}{4}\biggr)
T_l\biggl(\frac{\alpha'+\alpha''}{4}\biggr) C_{-\delta\epsilon}\biggl(\frac{\phi_1+\phi_2}{2}\biggr) \biggr]=0.
\end{gather*}
Similar computations, comparing equations \eqref{eq:segv1-v2} and \eqref{eq:primv2*-v1*}
show that
\begin{gather*}
 \biggl[T_l\biggl(\frac{\alpha^*-\alpha}{4}\biggr) S_{-\delta\epsilon}\biggl(\frac{\phi_2-\phi_1}{2}\biggr)-\delta\tau
T_l\biggl(\frac{\alpha'-\alpha''}{4}\biggr) S_{-\delta\epsilon}\biggl(\frac{\phi_2+\phi_1}{2}\biggr)\biggr] \\
 \qquad{}\times\biggl[ C_{-\delta\epsilon}\biggl(\frac{\phi_1+\phi_2}{2}\biggr)-\delta\tau l T_l\biggl(\frac{\alpha^*+\alpha}{4}\biggr)
T_l\biggl(\frac{\alpha'+\alpha''}{4}\biggr) C_{-\delta\epsilon}\biggl(\frac{\phi_1-\phi_2}{2}\biggr) \biggr]=0.
\end{gather*}
From the last two equations, we conclude that
\[
T_l\biggl(\frac{\alpha^*-\alpha}{4}\biggr) S_{-\delta\epsilon}\biggl(\frac{\phi_2-\phi_1}{2}\biggr)-\delta\tau
T_l\biggl(\frac{\alpha'-\alpha''}{4}\biggr) S_{-\delta\epsilon}\biggl(\frac{\phi_2+\phi_1}{2}\biggr)=0,
\]
which is exactly the superposition formula given by \eqref{eq:superpos}.

It remains to prove that $\alpha^*$ defined by \eqref{eq:superpos} satisfies \eqref{eqorig1_4}. We are assuming that
$\alpha'$ and $\alpha''$ satisfy \eqref{eq:btanalytic_v1} and \eqref{eq:btanalytic_v2}, respectively and $\alpha$ satisfies \eqref{eqorig1_4}. Theorem \ref{thBT} implies that $\alpha'$ and $\alpha''$ also satisfy \eqref{eqorig1_4}. We introduce the following notation:
\begin{gather}\label{Sigmaetall}
Z:= \frac{S_{-\delta\epsilon}\bigl(\frac{\phi_2+\phi_1}{2}\bigr) } {S_{-\delta\epsilon}\bigl(\frac{\phi_2-\phi_1}{2}\bigr)},\qquad D:=\alpha''-\alpha', \qquad \Sigma:=\alpha''+\alpha',
\qquad G:= -\delta\tau Z T_l \biggl(\frac{D}{4}\biggr).
\end{gather}
It follows from \eqref{eq:superpos} that
\begin{equation*}
\alpha^*_{x_ix_i}=\alpha_{x_ix_i}+\frac{4}{1+l G^2}\biggl(\frac{-2l G}{1+l G^2}(G_{x_i})^2+G_{x_ix_i}\biggr), \qquad i=1,2.
\end{equation*}
Therefore, proving that $\alpha^*$ satisfies \eqref{eqorig1_4} is equivalent to
proving
\begin{equation}\label{eqalst}
\frac{4}{1+l G^2}
\biggl(\frac{-2l G}{1+l G^2}\bigl[(G_{x_1})^2-(G_{x_2})^2\bigr]+G_{x_1x_1}-G_{x_2x_2} \biggr) +\epsilon\delta ( S_l(\alpha^*)-S_l(\alpha))=0,
\end{equation}
where we used the fact that $\alpha$ is a solution of \eqref{eqorig1_4}.

Observe that
\[
S_l(\alpha^*)-S_l(\alpha)=2S_l\biggl(\frac{\alpha^*-\alpha}{2}\biggr) C_l\biggl(\frac{\alpha^*+\alpha}{2}\biggr)=
\frac{4G}{1+l G^2}C_l\biggl(\frac{\alpha^*+\alpha}{2}\biggr),
\]
and it follows from \eqref{eq:superpos} and a long computation that
\[
C_l\biggl(\frac{\alpha^*+\alpha}{2}\biggr)=\frac{1}{1+l G^2}\bigl[ \bigl(1-l G^2\bigr)C_l(\alpha)-2l G S_l(\alpha)\bigr].
\]
Therefore, \eqref{eqalst} reduces to
\begin{gather}
-2l\bigl[(G_{x_1})^2-(G_{x_2})^2\bigr]+\frac{\bigl(1+l G^2\bigr)}{G}(G_{x_1x_1}-G_{x_2x_2}) \nonumber\\
\qquad{}+\epsilon\delta\bigl[ \bigl(1-l G^2\bigr)C_l(\alpha)-2l G S_l(\alpha)\bigr]=0.\label{eqalst1}
\end{gather}

Using the expression of $G$, we have
\[
G_{x_i}=\frac{G D_{x_i}}{2S_l\bigl(\frac{D}{2}\bigr)},
\quad
G_{x_ix_i}=\frac{G}{2S_l\bigl(\frac{D}{2}\bigr)}\biggl( D_{x_ix_i}+\frac{l}{2}T_l\biggl(\frac{D}{4}\biggr) (D_{x_i})^2 \biggr), \qquad i=1,2.
\]
Moreover,
\[
\frac{l}{4S_l\bigl(\frac{D}{2}\bigr)}\biggl(\frac{-2G^2}{S_l\bigl(\frac{D}{2}\bigr)}+\bigl(1+l G^2\bigr)
T_l\biggl(\frac{D}{4}\biggr)\biggr)= \frac{-l\bigl(Z^2-1\bigr)}{8\bigl(C_l\bigl(\frac{D}{4}\bigr)\bigr)^2}.
\]
Substituting these expressions into \eqref{eqalst1}, we conclude that we need to prove that
\begin{gather}
\frac{-l\bigl(Z^2-1\bigr)}{8C_l^2\bigl(\frac{D}{4}\bigr)} \bigl[(D_{x_1})^2-(D_{x_2})^2\bigr]+
\frac{1+l G^2}{2S_l\bigl(\frac{D}{2}\bigr)}(D_{x_1x_1}-D_{x_2x_2})\nonumber\\
\qquad{}+\epsilon\delta \bigl[ \bigl(1-l G^2\bigr)C_l(\alpha)-2l G S_l(\alpha)\bigr]=0.\label{eqalst3}
\end{gather}

Since $\alpha'$ and $\alpha''$ satisfy \eqref{eqorig1_4}, it follows that
\begin{equation}\label{eqdifDxixi}
D_{x_1x_1}-D_{x_2x_2}=-2\epsilon\delta S_l\biggl(\frac{D}{2}\biggr) C_l\biggl(\frac{\Sigma}{2}\biggr),
\end{equation}
where $\Sigma$ is given by \eqref{Sigmaetall}. Moreover, since $\alpha'$ and $\alpha''$ satisfy
\eqref{eq:btanalytic_v1} and \eqref{eq:btanalytic_v2} respectively,
 a long but straightforward computation shows that
\begin{gather*}
(D_{x_1})^2-(D_{x_2})^2=\frac{-4}{\lambda_1^2\lambda_2^2}
\biggl[ -\biggl( C_l^2\Bigl(\frac{\alpha}{2}\Bigr)S_l^2\biggl(\frac{\alpha'}{2}\biggr)-
C_l^2\biggl(\frac{\alpha'}{2}\biggr)S_l^2\Bigl(\frac{\alpha}{2}\Bigr)\biggr) (1-\Lambda_1^2)\lambda_2^2 \\
\qquad\qquad{}-\biggl( C_l^2\Bigl(\frac{\alpha}{2}\Bigr)S_l^2\biggl(\frac{\alpha''}{2}\biggr)-
C_l^2\biggl(\frac{\alpha''}{2}\biggr)S_l^2\Bigl(\frac{\alpha}{2}\Bigr)\biggr) (1-\Lambda_2^2)\lambda_1^2 \\
\qquad\qquad{}+2\biggl( C_l^2\Bigl(\frac{\alpha}{2}\Bigr)S_l\biggl(\frac{\alpha'}{2}\biggr)S_l\biggl(\frac{\alpha''}{2}\biggr) - C_l\biggl(\frac{\alpha'}{2}\biggr)C_l\biggl(\frac{\alpha''}{2}\biggr)S_l^2\Bigl(\frac{\alpha}{2}\Bigr) \biggr)(1-\Lambda_1\Lambda_2)\lambda_1\lambda_2 \\
\qquad\qquad{} -S_l(\alpha)S_l\biggl(\frac{D}{2}\biggr) (\Lambda_2-\Lambda_1)\delta\tau\lambda_1\lambda_2 \biggr].
\end{gather*}
Therefore, it follows from the fact
\[
1-\Lambda_i^2=-\epsilon\delta\lambda_i^2, \qquad
1-\Lambda_1\Lambda_2=\frac{-\epsilon\delta\bigl(Z^2+1\bigr)}{Z^2-1} \lambda_1\lambda_2, \qquad
\Lambda_2-\Lambda_1=\frac{2\epsilon\delta Z}{Z^2-1}\lambda_1\lambda_2,
\]
 that
\begin{gather*}
(D_{x_1})^2-(D_{x_2})^2=\frac{-16\epsilon\delta}{Z^2-1}\biggl[ \biggl( C_l^2\Bigl(\frac{\alpha}{2} \Bigr) C_l^2\biggl(\frac{\Sigma}{4}\biggr)- S_l^2\Bigl(\frac{\alpha}{2}\Bigr) S_l^2\biggl(\frac{\Sigma}{4}\biggr) \biggr)
S_l^2\biggl(\frac{D}{4}\biggr)Z^2 \\
\hphantom{(D_{x_1})^2-(D_{x_2})^2=}{}
-\biggl(\! C_l^2\Bigl(\frac{\alpha}{2}\Bigr)S_l^2\biggl(\frac{\Sigma}{4}\biggr)-
S_l^2\Bigl(\frac{\alpha}{2} \Bigr) C_l^2\biggl(\frac{\Sigma}{4}\biggr)\!\biggr)
C_l^2\biggl(\!\frac{D}{4}\!\biggr)
 -\frac{1}{2}S_l(\alpha)S_l\biggl(\!\frac{D}{2}\!\biggr)\delta\tau Z
 \biggr],
\end{gather*}
where we used the fact that $S_l^2\bigl(\frac{D}{4}\bigr)Z^2=G^2 C_l^2\bigl(\frac{D}{4}\bigr)$.

Substituting this expression and \eqref{eqdifDxixi} into \eqref{eqalst3}, we obtain
a second degree polynomial in~$G$, whose coefficients vanish identically.
\end{proof}

Theorem~\ref{th:superth1_4} shows that the composition of B\"acklund transformations, with distinct parameters, is commutative. Moreover, starting with a solution of the sine-Gordon equation
(resp.\ sinh-Gordon equation) on gets infinitely many solutions of the same equation, by an algebraic expression, after the fist integration of the B\"acklund transformation, as one can see in the following diagram:
\begin{equation}\label{diagram}
\begin{array}{@{}ccccccc}
 & & \alpha'& \stackrel{\phi_2}{\longrightarrow}& \bar{\alpha} & & \longrightarrow\\
 & \stackrel{\phi_1}{\nearrow}& & \stackrel{\phi_1}{\nearrow}& &\stackrel{\phi_3}{\searrow}\\
\alpha& \stackrel{\phi_2}{\longrightarrow}& \alpha''& & & & \hat{\alpha} \\
 & \stackrel{\phi_3}{\searrow}& & \stackrel{\phi_3}{\searrow} & &\stackrel{\phi_1}{\nearrow} \\
 & & \alpha'''& \stackrel{\phi_2}{\longrightarrow}& \tilde{\alpha} & &\longrightarrow
\end{array}
\end{equation}

Our next two theorems provide the superposition formulae for the elliptic sinh-Gordon and the elliptic sine-Gordon equations. The B\"acklund transformations for these equations, given in Corollary \ref{thBT5_6}, provide a solution of the elliptic sine-Gordon equation starting from a solution of the elliptic sinh-Gordon equation when $r=0$ and vice-versa when $r=1$. However, the superposition formula gives a new solution of the equation we started with, before applying any B\"acklund transformation.

\begin{Theorem}\label{th:super5}
Let $\alpha(x_1,x_2)$ be a solution of
\begin{equation}\label{elipsinh}
\alpha_{x_1x_1}+\alpha_{x_2x_2}=\sinh \alpha.
\end{equation}
Let $\alpha'$ and $\alpha''$ be solutions of the elliptic sine-Gordon equation associated to $\alpha$ by $\phi_1($resp.\ $\phi_2)$-B\"acklund transformation \eqref{eq:btanalytic_v5_6} where $r=0$ and $\tau=\pm 1$.
If $\phi_1\neq \phi_2$, then there exists a unique solution $\alpha^*$ of \eqref{elipsinh} which is associated to
$\alpha'$ and $\alpha''$ by $\phi_2($resp.\ $\phi_1)$-B\"acklund transformation~\eqref{eq:btanalytic_v5_6}, with $r=0$. Moreover, $\alpha^*$ is
given algebraically by
\begin{equation}\label{eq:superposcas5}
\tanh \frac{\alpha^*}{2}=
 \frac{P-Q \tanh \bigl({\frac{\alpha}{2}}\bigr)} {Q-P \tanh \bigl({\frac{\alpha}{2}}\bigr)},
\end{equation}
where
\begin{gather*}
P=\tau(\sinh\phi_2-\sinh\phi_1) \sin\biggl(\frac{\alpha''-\alpha'}{2}\biggr), \\
Q=(1+\sinh\phi_1\sinh\phi_2)\cos\biggl(\frac{\alpha''-\alpha'}{2}\biggr) -\cosh\phi_1\cosh\phi_2.
\end{gather*}
\end{Theorem}

\begin{proof} Let $\alpha$ be a solution of the elliptic sinh-Gordon equation, i.e., \eqref{eqorig5_6} with $r=0$.
Let $\alpha'$ and $\alpha''$ be solutions of the elliptic sine-Gordon equation, i.e., \eqref{eqprime5_6} with $r=0$, associated to $\alpha$
by $\phi_1$(resp.\ $\phi_2$)-B\"acklund transformation \eqref{eq:btanalytic_v5_6} with $r=0$, $\phi_1\neq \phi_2$. Then
 \begin{gather}
\alpha'_{x_1}-\alpha_{x_2}=\frac{2}{\cosh\phi_1}\sin \biggl(\frac{\alpha'}{2}\biggr)\cosh \Bigl(\frac{\alpha}{2}\Bigr)+\frac{2\tau\sinh\phi_1}{\cosh\phi_1} \cos\biggl(\frac{\alpha'}{2}\biggr)\sinh \Bigl(\frac{\alpha}{2}\Bigr),\nonumber\\
\alpha'_{x_2}+\alpha_{x_1}=-\frac{2}{\cosh\phi_1}\cos\biggl(\frac{\alpha'}{2}\biggr)\sinh\Bigl(\frac{\alpha}{2}\Bigr)+\frac{2\tau\sinh\phi_1}{\cosh\phi_1} \sin\biggl(\frac{\alpha'}{2}\biggr)\cosh\Bigl(\frac{\alpha}{2}\Bigr),\label{BTcase51}\\
\alpha''_{x_1}-\alpha_{x_2}=\frac{2}{\cosh\phi_2}\sin\biggl(\frac{\alpha''}{2}\biggr)\cosh \Bigl(\frac{\alpha}{2}\Bigr)+\frac{2\tau\sinh\phi_2}{\cosh\phi_2} \cos\biggl(\frac{\alpha''}{2}\biggr)\sinh\Bigl(\frac{\alpha}{2}\Bigr),\nonumber\\
\alpha''_{x_2}+\alpha_{x_1}=-\frac{2}{\cosh\phi_2}\cos\biggl(\frac{\alpha''}{2}\biggr)\sinh\Bigl(\frac{\alpha}{2}\Bigr)+\frac{2\tau\sinh\phi_2}{\cosh\phi_2} \sin\biggl(\frac{\alpha''}{2}\biggr)\cosh\Bigl(\frac{\alpha}{2}\Bigr).\label{BTcase52}
\end{gather}
Suppose $\alpha^*$ exists then it should satisfy the following
systems of equations
\begin{gather}
\alpha^*_{x_1}-\alpha'_{x_2}=\frac{2}{\cosh\phi_2}\sinh\biggl(\frac{\alpha^*}{2}\biggr)\cos\biggl(\frac{\alpha'}{2}\biggr)+\frac{2\tau\sinh\phi_2}{\cosh\phi_2} \cosh\biggl(\frac{\alpha^*}{2}\biggr)\sin\biggl(\frac{\alpha'}{2}\biggr),\nonumber\\
\alpha^*_{x_2}+\alpha'_{x_1}=-\frac{2}{\cosh\phi_ 2}\cosh\biggl(\frac{\alpha^*}{2}\biggr)\sin\biggl(\frac{\alpha'}{2}\biggr)+\frac{2\tau\sinh\phi_2}{\cosh\phi_2} \sinh\biggl(\frac{\alpha^*}{2}\biggr)\cos\biggl(\frac{\alpha'}{2}\biggr),\label{BTcase62}
\\
\alpha^*_{x_1}-\alpha''_{x_2}=\frac{2}{\cosh\phi_1}\sinh\biggl(\frac{\alpha^*}{2}\biggr)\cos\biggl(\frac{\alpha''}{2}\biggr)+\frac{2\tau\sinh\phi_1}{\cosh\phi_1} \cosh\biggl(\frac{\alpha^*}{2}\biggr)\sin\biggl(\frac{\alpha''}{2}\biggr),\nonumber\\
\alpha^*_{x_2}+\alpha''_{x_1}=-\frac{2}{\cosh\phi_1}\cosh\biggl(\frac{\alpha^*}{2}\biggr)\sin\biggl(\frac{\alpha''}{2}\biggr)+\frac{2\tau\sinh\phi_1}{\cosh\phi_1} \sinh\biggl(\frac{\alpha^*}{2}\biggr)\cos\biggl(\frac{\alpha''}{2}\biggr).\label{BTcase61}
\end{gather}
Subtracting the first (resp.\ second) equation of \eqref{BTcase52} from the first (resp.\ second) equation of~\eqref{BTcase51}, we get the system of equations
\begin{gather}
\alpha'_{x_1}-\alpha''_{x_1}=\frac{2}{\cosh\phi_1}\sin\biggl(\frac{\alpha'}{2}\biggr)\cosh \Bigl(\frac{\alpha}{2}\Bigr)+\frac{2\tau\sinh\phi_1}{\cosh\phi_1} \cos\biggl(\frac{\alpha'}{2}\biggr)\sinh\Bigl(\frac{\alpha}{2}\Bigr)
\nonumber\\ \hphantom{\alpha'_{x_1}-\alpha''_{x_1}=}{}
- \frac{2}{\cosh\phi_2}\sin\biggl(\frac{\alpha''}{2}\biggr)\cosh \Bigl(\frac{\alpha}{2}\Bigr)-\frac{2\tau\sinh\phi_2}{\cosh\phi_2} \cos\biggl(\frac{\alpha''}{2}\biggr)\sinh\Bigl(\frac{\alpha}{2}\Bigr),\nonumber\\
\alpha'_{x_2}-\alpha''_{x_2}= -\frac{2}{\cosh\phi_1}\cos\biggl(\frac{\alpha'}{2}\biggr)\sinh\Bigl(\frac{\alpha}{2}\Bigr)+\frac{2\tau\sinh\phi_1}{\cosh\phi_1} \sin\biggl(\frac{\alpha'}{2}\biggr)\cosh\Bigl(\frac{\alpha}{2}\Bigr)
\nonumber\\ \hphantom{\alpha'_{x_2}-\alpha''_{x_2}=}{}
+\frac{2}{\cosh\phi_2}\cos\biggl(\frac{\alpha''}{2}\biggr)\sinh\Bigl(\frac{\alpha}{2}\Bigr)-\frac{2\tau\sinh\phi_2}{\cosh\phi_2} \sin\biggl(\frac{\alpha''}{2}\biggr)\cosh\Bigl(\frac{\alpha}{2}\Bigr). \label{BT51minusBT52}
\end{gather}

Similarly, subtracting the first (resp.\ second) equation of \eqref{BTcase61} from the first (resp.\ second) equation of \eqref{BTcase62}, we get
\begin{gather}
\alpha'_{x_2}-\alpha''_{x_2}=-\frac{2}{\cosh\phi_2}\sinh\biggl(\frac{\alpha^*}{2}\biggr)\cos\biggl(\frac{\alpha'}{2}\biggr)-\frac{2\tau\sinh\phi_2}{\cosh\phi_2} \cosh\biggl(\frac{\alpha^*}{2}\biggr)\sin\biggl(\frac{\alpha'}{2}\biggr)
\nonumber\\ \hphantom{\alpha'_{x_2}-\alpha''_{x_2}=}{}
 +\frac{2}{\cosh\phi_1}\sinh\biggl(\frac{\alpha^*}{2}\biggr)\cos\biggl(\frac{\alpha''}{2}\biggr)+\frac{2\tau\sinh\phi_1}{\cosh\phi_1} \cosh\biggl(\frac{\alpha^*}{2}\biggr)\sin\biggl(\frac{\alpha''}{2}\biggr),\nonumber\\
\alpha'_{x_1}-\alpha''_{x_1}=-\frac{2}{\cosh\phi_ 2}\cosh\biggl(\frac{\alpha^*}{2}\biggr)\sin\biggl(\frac{\alpha'}{2}\biggr)+\frac{2\tau\sinh\phi_2}{\cosh\phi_2} \sinh\biggl(\frac{\alpha^*}{2}\biggr)\cos\biggl(\frac{\alpha'}{2}\biggr)
\nonumber\\ \hphantom{\alpha'_{x_1}-\alpha''_{x_1}=}{}
 +\frac{2}{\cosh\phi_1}\cosh\biggl(\frac{\alpha^*}{2}\biggr)\sin\biggl(\frac{\alpha''}{2}\biggr)- \frac{2\tau\sinh\phi_1}{\cosh\phi_1} \sinh\biggl(\frac{\alpha^*}{2}\biggr)\cos\biggl(\frac{\alpha''}{2}\biggr). \label{BT62minusBT61}
\end{gather}

Comparing equation \eqref{BT51minusBT52} with \eqref{BT62minusBT61}, we obtain two algebraic
expressions that must be satisfied by $\cosh(\alpha^*/2)$ and $\sinh(\alpha^*/2)$.
By solving this system of algebraic equations, we conclude that $\alpha^*/2$ is given by
\eqref{eq:superposcas5}.

Now we need to prove that $\alpha^*$ satisfies the elliptic sinh-Gordon equation \eqref{elipsinh}.
Recall that~$\alpha'$ and~$\alpha''$ are solutions of the elliptic sine-Gordon equation.
We introduce the following notation
\begin{gather}
A=\tau(\sinh\phi_2-\sinh\phi_1),\qquad B=\cosh\phi_1\cosh\phi_2,\qquad L=1+\sinh\phi_1\sinh\phi_2, \nonumber\\
 D=\alpha''-\alpha'.\label{ABLD}
\end{gather}
Then the functions $P$ and $Q$ are rewritten as
\begin{equation*}
P=A\sin\biggl(\frac{D}{2}\biggr),\qquad Q=L\cos\biggl(\frac{D}{2}\biggr)-B \qquad \mbox {and}\qquad B^2-L^2=A^2.
\end{equation*}
We denote the right-hand side of \eqref{eq:superposcas5} by
\begin{equation*}
H=\frac{G-\tanh\bigl(\frac{\alpha}{2}\bigr)}{1-G\tanh\bigl(\frac{\alpha}{2}\bigr)}, \qquad \mbox{where}\quad G=\frac{P}{Q}.
\end{equation*}
It follows that proving that $\alpha^*$ satisfies \eqref{elipsinh} is equivalent to proving that
\begin{equation}\label{eqH}
\frac{2H}{1-H^2}|\nabla H|^2+ \triangle H-H=0.
\end{equation}
Since $\alpha$ satisfies \eqref{elipsinh} it follows that
\begin{gather*}
\triangle \tanh \Bigl(\frac{\alpha}{2}\Bigr)= \frac{1}{2\cosh^2\bigl(\frac{\alpha}{2}\bigr)}\Bigl(\sinh \alpha-\tanh\Bigl(\frac{\alpha}{2}\Bigr)|\nabla\alpha|^2\Bigr), \\
\Bigl|\nabla\tanh\Bigl(\frac{\alpha}{2}\Bigr)\Bigr|^2=\frac{1}{4\cosh^4 \bigl(\frac{\alpha}{2}\bigr)}
|\nabla \alpha|^2.
\end{gather*}
Therefore, proving \eqref{eqH} reduces to showing
that $G$ satisfies
\begin{equation}\label{eqG}
\frac{2G}{1-G^2}|\nabla G|^2+\triangle G-G(\cosh \alpha-G\sinh\alpha)=0.
\end{equation}

Long but straightforward computations show that
\begin{gather*}
 \triangle G = \frac{A}{2 Q^2}\biggl(L-B\cos\biggl(\frac{D}{2}\biggr) \biggr) {\triangle D}
-\frac{A}{4Q^3}\bigl(BQ+2A^2\bigr)\sin \biggl(\frac{D}{2}\biggr) |\nabla D|^2, \\
 |\nabla G|^2 = \frac{A^2}{4 Q^4}\biggl(L-B\cos\biggl(\frac{D}{2}\biggr) \biggr)^2 |\nabla D|^2,\\
 1-G^2 = \frac{1}{Q^2}\biggl(L-B\cos\biggl(\frac{D}{2}\biggr) \biggr)^2.
\end{gather*}
Moreover, since $D=\alpha''-\alpha'$ where $\alpha'$ and $\alpha''$ are solutions of the elliptic sine-Gordon equation, it follows that
\begin{eqnarray*}
 \triangle D=\sin \alpha''-\sin \alpha'= 2\sin\frac{D}{2}\cos\biggl(\frac{\alpha''+\alpha'}{2}\biggr).
\end{eqnarray*}
Hence proving \eqref{eqG} reduces to proving
\begin{equation}\label{eqD}
-B|\nabla D|^2+\frac{2}{\sin\bigl(\frac{D}{2}\bigr)}\biggl(L-B\cos\biggl(\frac{D}{2}\biggr)\biggr)\triangle D-4Q\cosh
\alpha+4A\sin\biggl(\frac{D}{2}\biggr)\sinh\alpha=0.
\end{equation}

Now $\alpha'$ and $\alpha''$ satisfy \eqref{BTcase51} and \eqref{BTcase52} respectively.
These systems of equations provide
\begin{gather*}
 |\nabla D|^2= \frac{4}{B}
 \biggl[\cosh \alpha \biggl( B-L\cos\biggl(\frac{D}{2}\biggr)\biggr) +
 \cos \biggl(\frac{\alpha''+\alpha'}{2}\biggr)\biggl( L-B\cos\biggl(\frac{D}{2}\biggr)\biggr)
 \nonumber \\ \hphantom{|\nabla D|^2= \frac{4}{B}\biggl[}{}
 + A\sin \biggl(\frac{D}{2}\biggr)\sinh \alpha\biggr].
\end{gather*}
This last expression proves that \eqref{eqD} is identically satisfied, which concludes the proof that
$\alpha^*$ satisfies the elliptic sinh-Gordon equation
\end{proof}

The next result provides the superposition formula for the B\"acklund transformation \eqref{eq:btanalytic_v5_6}, when $r=1$, which gives solutions for the elliptic sine-Gordon equation.

\begin{Theorem}\label{th:super6}
Let $\alpha(x_1,x_2)$ be a solution of
\begin{equation}\label{elipsin_th3}
\alpha_{x_1x_1}+\alpha_{x_2x_2}=\sin \alpha.
\end{equation}
Let $\alpha'$ and $\alpha''$ be solutions of the elliptic sinh-Gordon
equation associated to $\alpha$ by a $\phi_1($resp.~$\phi_2)$-B\"acklund transformation \eqref{eq:btanalytic_v5_6}, where $r=1$ and $\tau=\pm 1$.
If $\phi_1\neq \phi_2$, then there exists a~unique solution $\alpha^*$ of \eqref{elipsin_th3} which is associated to
$\alpha'$ and $\alpha''$ by the $\phi_2($resp.\ $\phi_1)$-B\"acklund transformation \eqref{eq:btanalytic_v5_6}, with $r=1$. Moreover, $\alpha^*$ is
given algebraically by
\begin{equation}\label{eq:superpos_th3}
\tan \frac{\alpha^*}{2}= \frac{\tilde{P}-\tilde{Q} {\tan \bigl(\frac{\alpha}{2}\bigr) } }{\tilde{Q}+\tilde{P} {\tan\bigl(\frac{\alpha}{2}\bigr)} },
\end{equation}
where
\begin{gather}\label{Pt}
\tilde{P}=-\tau(\sinh \phi_1-\sinh \phi_2) \sinh\biggl(\frac{\alpha''-\alpha'}{2}\biggr),
\\
\label{Qt}
\tilde{Q}=(1+\sinh\phi_1\sinh\phi_2)\cosh\frac{\alpha''-\alpha'}{2}-\cosh\phi_1\cosh\phi_2.
\end{gather}
\end{Theorem}

\begin{proof} The proof is similar to the previous one.
Let $\alpha$ be a solution of the elliptic sine-Gordon equation, i.e., \eqref{eqorig5_6} with $r=1$.
Let $\alpha'$ and $\alpha''$ be solutions of the
elliptic sinh-Gordon (i.e.,~\eqref{eqprime5_6} with $r=1$)
associated to $\alpha$ by $\phi_1$(resp.\ $\phi_2$)-B\"acklund transformation
\eqref{eq:btanalytic_v5_6} with $r=1$, where $\phi_1\neq \phi_2$, then
\begin{gather}
\alpha'_{x_1}-\alpha_{x_2}=\frac{2}{\cosh\phi_1}\sinh \biggl(\frac{\alpha'}{2}\biggr)\cos \Bigl(\frac{\alpha}{2}\Bigr)+\frac{2\tau\sinh\phi_1}{\cosh\phi_1} \cosh\biggl(\frac{\alpha'}{2}\biggr)\sin \Bigl(\frac{\alpha}{2}\Bigr),\nonumber\\
\alpha'_{x_2}+\alpha_{x_1}=-\frac{2}{\cosh\phi_1}\cosh\biggl(\frac{\alpha'}{2}\biggr)\sin\Bigl(\frac{\alpha}{2}\Bigr)
+\frac{2\tau\sinh\phi_1}{\cosh\phi_1} \sinh\biggl(\frac{\alpha'}{2}\biggr)\cos\Bigl(\frac{\alpha}{2}\Bigr),\label{BTcase61_th3}\\
\alpha''_{x_1}-\alpha_{x_2}=\frac{2}{\cosh\phi_2}\sinh\biggl(\frac{\alpha''}{2}\biggr)\cos \Bigl(\frac{\alpha}{2}\Bigr)
+\frac{2\tau\sinh\phi_2}{\cosh\phi_2} \cosh\biggl(\frac{\alpha''}{2}\biggr)\sin\Bigl(\frac{\alpha}{2}\Bigr),\nonumber\\
\alpha''_{x_2}+\alpha_{x_1}=-\frac{2}{\cosh\phi_2}\cosh\biggl(\frac{\alpha''}{2}\biggr)\sin\Bigl(\frac{\alpha}{2}\Bigr)
+\frac{2\tau\sinh\phi_2}{\cosh\phi_2} \sinh\biggl(\frac{\alpha''}{2}\biggr)\cos\Bigl(\frac{\alpha}{2}\Bigr).\label{BTcase62_th3}
\end{gather}
Suppose $\alpha^*$ exists then it should satisfy the following
systems of equations:
\begin{gather}
\alpha^*_{x_1}-\alpha'_{x_2}=\frac{2}{\cosh\phi_2}\sin\biggl(\frac{\alpha^*}{2}\biggr)\cosh\biggl(\frac{\alpha'}{2}\biggr)+\frac{2\tau\sinh\phi_2}{\cosh\phi_2} \cos\biggl(\frac{\alpha^*}{2}\biggr)\sinh\biggl(\frac{\alpha'}{2}\biggr),\nonumber\\
\alpha^*_{x_2}+\alpha'_{x_1}=-\frac{2}{\cosh\phi_ 2}\cos\biggl(\frac{\alpha^*}{2}\biggr)\sinh\biggl(\frac{\alpha'}{2}\biggr)+\frac{2\tau\sinh\phi_2}{\cosh\phi_2} \sin\biggl(\frac{\alpha^*}{2}\biggr)\cosh\biggl(\frac{\alpha'}{2}\biggr),\label{BTcase52_th3}\\
\alpha^*_{x_1}-\alpha''_{x_2}=\frac{2}{\cosh\phi_1}\sin\biggl(\frac{\alpha^*}{2}\biggr)\cosh\biggl(\frac{\alpha''}{2}\biggr)+\frac{2\tau\sinh\phi_1}{\cosh\phi_1} \cos\biggl(\frac{\alpha^*}{2}\biggr)\sinh\biggl(\frac{\alpha''}{2}\biggr),\nonumber\\
\alpha^*_{x_2}+\alpha''_{x_1}=-\frac{2}{\cosh\phi_1}\cos\biggl(\frac{\alpha^*}{2}\biggr)\sinh\biggl(\frac{\alpha''}{2}\biggr)+\frac{2\tau\sinh\phi_1}{\cosh\phi_1} \sin\biggl(\frac{\alpha^*}{2}\biggr)\cosh\biggl(\frac{\alpha''}{2}\biggr).\label{BTcase51_th3}
\end{gather}
Subtracting the first (resp.\ second) equation of \eqref{BTcase52_th3} from the first (resp.\ second) equation of~\eqref{BTcase51_th3}, we get the system of equations
\begin{gather}
\alpha'_{x_1}-\alpha''_{x_1}=\frac{2}{\cosh\phi_1}\sinh\biggl(\frac{\alpha'}{2}\biggr)\cos \Bigl(\frac{\alpha}{2}\Bigr)+\frac{2\tau\sinh\phi_1}{\cosh\phi_1} \cosh\biggl(\frac{\alpha'}{2}\biggr)\sin\Bigl(\frac{\alpha}{2}\Bigr)
\nonumber\\ \hphantom{\alpha'_{x_1}-\alpha''_{x_1}=}{}
 - \frac{2}{\cosh\phi_2}\sinh\biggl(\frac{\alpha''}{2}\biggr)\cos \Bigl(\frac{\alpha}{2}\Bigr)-\frac{2\tau\sinh\phi_2}{\cosh\phi_2} \cosh\biggl(\frac{\alpha''}{2}\biggr)\sin\Bigl(\frac{\alpha}{2}\Bigr),\nonumber\\
\alpha'_{x_2}-\alpha''_{x_2}= -\frac{2}{\cosh\phi_1}\cosh\biggl(\frac{\alpha'}{2}\biggr)\sin\Bigl(\frac{\alpha}{2}\Bigr)+\frac{2\tau\sinh\phi_1}{\cosh\phi_1} \sinh\biggl(\frac{\alpha'}{2}\biggr)\cos\Bigl(\frac{\alpha}{2}\Bigr)
\nonumber\\ \hphantom{\alpha'_{x_2}-\alpha''_{x_2}=}{}
+\frac{2}{\cosh\phi_2}\cosh\biggl(\frac{\alpha''}{2}\biggr)\sin\Bigl(\frac{\alpha}{2}\Bigr)-\frac{2\tau\sinh\phi_2}{\cosh\phi_2} \sinh\biggl(\frac{\alpha''}{2}\biggr)\cos\Bigl(\frac{\alpha}{2}\Bigr).\label{BT61minusBT62th3}
\end{gather}

Similarly, subtracting the first (resp.\ second) equation of \eqref{BTcase51_th3} from the first (resp.\ second) equation of \eqref{BTcase52_th3}, we get
\begin{gather} 
\alpha'_{x_2}-\alpha''_{x_2}=-\frac{2}{\cosh\phi_2}\sin\biggl(\frac{\alpha^*}{2}\biggr)\cosh\biggl(\frac{\alpha'}{2}\biggr)-\frac{2\tau\sinh\phi_2}{\cosh\phi_2} \cos\biggl(\frac{\alpha^*}{2}\biggr)\sinh\biggl(\frac{\alpha'}{2}\biggr)
\nonumber\\ \hphantom{\alpha'_{x_2}-\alpha''_{x_2}=}{}
 +\frac{2}{\cosh\phi_1}\sin\biggl(\frac{\alpha^*}{2}\biggr)\cosh\biggl(\frac{\alpha''}{2}\biggr)+\frac{2\tau\sinh\phi_1}{\cosh\phi_1} \cos\biggl(\frac{\alpha^*}{2}\biggr)\sinh\biggl(\frac{\alpha''}{2}\biggr),\nonumber\\
\alpha'_{x_1}-\alpha''_{x_1}=-\frac{2}{\cosh\phi_ 2}\cos\biggl(\frac{\alpha^*}{2}\biggr)\sinh\biggl(\frac{\alpha'}{2}\biggr)+\frac{2\tau\sinh\phi_2}{\cosh\phi_2} \sin\biggl(\frac{\alpha^*}{2}\biggr)\cosh\biggl(\frac{\alpha'}{2}\biggr)
\nonumber\\ \hphantom{\alpha'_{x_1}-\alpha''_{x_1}=}{}
 +\frac{2}{\cosh\phi_1}\cos\biggl(\frac{\alpha^*}{2}\biggr)\sinh\biggl(\frac{\alpha''}{2}\biggr)- \frac{2\tau\sinh\phi_1}{\cosh\phi_1} \sin\biggl(\frac{\alpha^*}{2}\biggr)\cosh\biggl(\frac{\alpha''}{2}\biggr).
\label{BT52minusBT51_th3}
\end{gather}

Comparing equation \eqref{BT61minusBT62th3} with \eqref{BT52minusBT51_th3}, we obtain two algebraic
expressions that must be satisfied by $\cos(\alpha^*/2)$ and $\sin(\alpha^*/2)$.
By solving this system of equations, we conclude that~$\alpha^*/2$ is given by \eqref{eq:superpos_th3}, where $\tilde{P}$ and $\tilde{Q}$ are given by \eqref{Pt} and \eqref{Qt}.

 We need to prove that $\alpha^*$ defined by \eqref{eq:superpos_th3} satisfies the elliptic sine-Gordon equation \eqref{elipsin_th3}.
Recall that $\alpha'$ and $\alpha''$ are solutions of the elliptic sinh-Gordon equation. We use the previous notation \eqref{ABLD} for the constants
$A$, $B$, $L$ and the function $D=\alpha''-\alpha'$.
Then the functions $\tilde{P}$ and $\tilde{Q}$ are rewritten as
\begin{equation*}
\tilde{P}=A\sinh\biggl(\frac{D}{2}\biggr),\qquad \tilde{Q}=L\cosh\biggl(\frac{D}{2}\biggr)-B
\end{equation*}
and $B^2-L^2=A^2$.
We denote the right-hand side of \eqref{eq:superpos_th3} by
\begin{equation*}
\tilde{H}=\frac{\tilde{G}-\tan\bigl(\frac{\alpha}{2}\bigr)}{1+\tilde{G}\tan\bigl(\frac{\alpha}{2}\bigr)}, \qquad \mbox{where}\quad \tilde{G}=\frac{\tilde{P}}{\tilde{Q}}.
\end{equation*}
It follows that proving that $\alpha^*$ satisfies \eqref{elipsin_th3} is equivalent to proving that
\begin{equation}\label{eqHt}
\frac{2\tilde{H}}{1+\tilde{H}^2}|\nabla \tilde{H}|^2- \triangle \tilde{H}+\tilde{H}=0.
\end{equation}

Since $\alpha$ satisfies \eqref{elipsin_th3}, it follows that
\begin{gather*}
\triangle \tan \Bigl(\frac{\alpha}{2}\Bigr)= \frac{1}{2\cos^3\bigl(\frac{\alpha}{2}\bigr)}\Bigl(\sin \Bigl(\frac{\alpha}{2}\Bigr) |\nabla\alpha|^2+\cos\Bigl(\frac{\alpha}{2}\Bigr)\sin\alpha\Bigr), \\
\Bigl|\nabla\tan\Bigl(\frac{\alpha}{2}\Bigr)\Bigr|^2=\frac{1}{4\cos^4 \bigl(\frac{\alpha}{2}\bigr)}
|\nabla \alpha|^2.
\end{gather*}
Therefore, proving \eqref{eqHt} reduces to showing that $\tilde{G}$
satisfies
\begin{equation}\label{eqGt}
\triangle \tilde{G}-\frac{2\tilde{G}}{1+\tilde{G}^2}|\nabla\tilde{G}|^2-\tilde{G}(\cos \alpha +\tilde{G}\sin \alpha)=0.
\end{equation}
A long but straightforward computation shows that
\begin{gather*}
 \triangle \tilde{G} = \frac{A}{2 \tilde{Q}^2}\biggl(L-B\cosh\biggl(\frac{D}{2}\biggr) \biggr) {\triangle D}
+\frac{A}{4\tilde{Q}^3}\bigl(B\tilde{Q}-2A^2\bigr)\sinh\biggl(\frac{D}{2}\biggr)|\nabla D|^2,\\
 |\nabla \tilde{G}|^2 = \frac{A^2}{4 \tilde{Q}^4}\biggl(L-B\cosh\biggl(\frac{D}{2}\biggr) \biggr)^2|\nabla D|^2,\\
 1+\tilde{G}^2 = \frac{1}{\tilde{Q}^2}\biggl(L-B\cosh\biggl(\frac{D}{2}\biggr) \biggr)^2.
\end{gather*}

Since $\alpha'$ and $\alpha''$ are solutions of the elliptic sinh-Gordon equation and $D=\alpha''-\alpha'$, it follows that
\begin{eqnarray*}
 \triangle D=\sinh \alpha''-\sinh \alpha'=2\sinh\biggl(\frac{D}{2}\biggr)\cosh\biggl(\frac{\alpha''+\alpha'}{2}\biggr).
\end{eqnarray*}
Therefore, proving \eqref{eqGt} reduces to proving
\begin{equation}\label{eqD_th3}
\frac{B}{4}|\nabla D|^2+\biggl(L-B\cosh\biggl(\frac{D}{2}\biggr)\biggr)\cosh\biggl(\frac{\alpha''+\alpha'}{2}\biggr)-\tilde{P}\sin \alpha-\tilde{Q}\cos\alpha=0.
\end{equation}

Moreover, $\alpha'$ and $\alpha''$ satisfy \eqref{BTcase61_th3} and \eqref{BTcase62_th3} respectively, which provide $\nabla \alpha'$ and
$\nabla \alpha''$ and~$|\nabla D|^2$. Finally,
substituting this expression into \eqref{eqD_th3}, we obtain an identity.
Therefore, we conclude that $\alpha^*$ satisfies the elliptic sine-Gordon equation.
\end{proof}

The composition of B\"acklund transformations for the elliptic sinh-Gordon equation (ESHG) and the elliptic sine-Gordon equation (ESGE) is illustrated in the following diagrams. Starting with a solution of one of these equations, we need the superposition formula (i.e., a composition of two B\"acklund transformations) to produce new solutions of the same equation:
\[
\begin{array}{@{}c@{\,}c@{\,}c@{\,}c@{\,}c@{\,}c@{\,}c@{}}
 {\color{red}ESHG} & & {\color{blue}ESG} & & {\color{red}ESHG} & & {\color{blue}ESG} \\
 & & \alpha_1& \stackrel{\phi_2}{\longrightarrow}& {\alpha_{12}} & & \longrightarrow\\
 & \stackrel{\phi_1}{\nearrow}& & \stackrel{\phi_1}{\nearrow}& &\stackrel{\phi_3}{\searrow}\\
 \alpha & \stackrel{\phi_2}{\longrightarrow} & \alpha_2& & & & \alpha_{123} \\
 & \stackrel{\phi_3}{\searrow}& & \stackrel{\phi_3}{\searrow} & &\stackrel{\phi_1}{\nearrow} \\
 & & \alpha_3& \stackrel{\phi_2}{\longrightarrow}& {\alpha_{23}} & & \longrightarrow
\end{array}
\qquad
\begin{array}{ccccccc}
{\color{blue}ESG} & & {\color{red}ESHG} & & {\color{blue}ESG} & & \! \! {\color{red}ESHG} \\
 & & \alpha_1& \stackrel{\phi_2}{\longrightarrow}& {\alpha_{12}} & & \longrightarrow\\
 & \stackrel{\phi_1}{\nearrow}& & \stackrel{\phi_1}{\nearrow}& &\stackrel{\phi_3}{\searrow}\\
\alpha& \stackrel{\phi_2}{\longrightarrow}& \alpha_{2}& & & & {\alpha_{123}} \\
 & \stackrel{\phi_3}{\searrow}& & \stackrel{\phi_3}{\searrow} & &\stackrel{\phi_1}{\nearrow} \\
 & & \alpha_3& \stackrel{\phi_2}{\longrightarrow}& \alpha_{23} & & \longrightarrow
\end{array}
\]

\section{Six cases}\label{sec:6cases} In this section, aiming the applications of Section \ref{examples} and also in order to help the reader to easily access the theory in each case, we give the explicit
B\"acklund transformation and the corresponding superposition formula, for each one of the six cases considered in the previous sections.

{\bf Case 1.}
Considering $s=0$, $\epsilon=1$, $r=0$ and $\delta=-1$ then $l=1$ and \eqref{eq:lambdaephi} implies that $\Lambda^2+\lambda^2=1$ hence,
$\Lambda=\cos \phi$ and $\lambda=\sin\phi$, where $\phi\in (0,\pi)$.

\begin{Theorem}\label{SG}
For any solution $\alpha(x_1,x_2)$ of
\begin{equation}\label{eq:SG}
\alpha_{x_1x_1}-\alpha_{x_2x_2}=\sin \alpha,
\end{equation}
 the following system of differential equations for $\alpha'(x_1,x_2)$, where $\phi\in (0,\pi)$, is integrable:
\begin{gather}
\alpha'_{x_1}+\alpha_{x_2}=\frac{2}{\sin\phi}\sin\biggl(\frac{\alpha'}{2}\biggr)\cos\Bigl(\frac{\alpha}{2}\Bigr)-\frac{2\tau \cos\phi}{\sin\phi} \cos\biggl(\frac{\alpha'}{2}\biggr)\sin\Bigl(\frac{\alpha}{2}\Bigr),\nonumber\\
\alpha'_{x_2}+\alpha_{x_1}=-\frac{2}{\sin\phi}\cos\biggl(\frac{\alpha'}{2}\biggr)\sin\Bigl(\frac{\alpha}{2}\Bigr)+\frac{2\tau \cos\phi}{\sin\phi} \sin\biggl(\frac{\alpha'}{2}\biggr)\cos\Bigl(\frac{\alpha}{2}\Bigr).\label{Backcase1}
\end{gather}
Moreover, $\alpha'$ is also a solution of \eqref{eq:SG}. If $\alpha'$ and $\alpha''$ are solutions of \eqref{eq:SG} obtained by the B\"acklund transformation \eqref{Backcase1}, for constants $\phi_1\neq \phi_2$, then a fourth solution $\alpha^*$ of \eqref{eq:SG} is obtained algebraically with
 the superposition formula
 \begin{equation}\label{superposcase1}
\tan\biggl(\frac{\alpha^*-\alpha}{4}\biggr) =
-\tau \frac{\sin\bigl(\frac{\phi_2+\phi_1}{2}\bigr) } {\sin \bigl(\frac{\phi_2-\phi_1}{2}\bigr)}
\tan\biggl(\frac{\alpha'-\alpha''}{4}\biggr).
 \end{equation}
\end{Theorem}

This theorem gives the classical B\"acklund transformation \eqref{Backcase1} and the superposition formula~\eqref{superposcase1} for the sine-Gordon equation obtained by B\"acklund \cite{Backlund1875,Backlund1905} and Bianchi \cite{bianchi1}, when $\tau=-1$. It can also be
found for example in \cite[p.~44]{Tenenblat1998}.

{\bf Case 2.}
Considering $s=1$, $\epsilon=1$, $r=0$ and $\delta=1$ then $l=1$ and \eqref{eq:lambdaephi} implies that
$\Lambda^2 -\lambda^2=1$, hence $\Lambda=\cosh\phi$ and $\lambda=\sinh\phi$,
where $\phi\in (0,\infty)$.

\begin{Theorem}
For any solution $\alpha(x_1,x_2)$ of
\begin{equation}\label{eq:SGmin}
\alpha_{x_1x_1}-\alpha_{x_2x_2}=-\sin \alpha,
\end{equation}
 the following system of differential equations for $\alpha'(x_1,x_2)$, $\phi\in(0,\infty)$, is integrable:
\begin{gather}
\alpha'_{x_1}+\alpha_{x_2}=\frac{2}{\sinh\phi}\sin\biggl(\frac{\alpha'}{2}\biggr)\cos\Bigl(\frac{\alpha}{2}\Bigr) +\frac{2\tau\cosh\phi}{\sinh\phi}\cos\biggl(\frac{\alpha'}{2}\biggr)\sin\Bigl(\frac{\alpha}{2}\Bigr),\nonumber\\
\alpha'_{x_2}+\alpha_{x_1}=-\frac{2}{\sinh\phi}\cos\biggl(\frac{\alpha'}{2}\biggr)\sin\Bigl(\frac{\alpha}{2}\Bigr) -\frac{2\tau\cosh\phi}{\sinh\phi}\sin\biggl(\frac{\alpha'}{2}\biggr)\cos\Bigl(\frac{\alpha}{2}\Bigr).\label{Backcase2}
\end{gather}
 Moreover, $\alpha'$ is also a solution of \eqref{eq:SGmin}.
If $\alpha'$ and $\alpha''$ are solutions of \eqref{eq:SGmin} obtained by the B\"acklund transformation \eqref{Backcase2}, for constants $\phi_1\neq \phi_2$, then a fourth solution $\alpha^*$ of \eqref{eq:SGmin} is obtained algebraically with the superposition formula
\begin{equation*}
\tan \biggl(\frac{\alpha^*-\alpha}{4}\biggr) =
\tau \frac{\sinh\bigl(\frac{\phi_2+\phi_1}{2}\bigr) } {\sinh\bigl(\frac{\phi_2-\phi_1}{2}\bigr)}
\tan\biggl(\frac{\alpha'-\alpha''}{4}\biggr).
\end{equation*}
\end{Theorem}

We observe that Cases 1 and 2 deal with the sine-Gordon equation and they are complementary in the following sense: while in Case 1 the parameter $\phi\in(0,\pi)$ and hence
$\Lambda=\cos\phi\in (-1,1)$, in Case 2, the parameter $\phi\in (0,\infty)$ and hence $\Lambda=\cosh \phi \in (1,+\infty)$.

{\bf Case 3.}
Considering $s=1$, $\epsilon=1$, $r=1$ and $\delta=1$ then $l=-1$ and \eqref{eq:lambdaephi} implies that
$\Lambda^2 -\lambda^2=1$, hence $\Lambda=\cosh\phi$ and $\lambda=\sinh\phi$,
where $\phi\in (0,\infty)$.

\begin{Theorem}
For any solution $\alpha(x_1,x_2)$ of
\begin{equation}\label{eq:SHGmin}
\alpha_{x_1x_1}-\alpha_{x_2x_2}=-\sinh \alpha,
\end{equation}
 the following system of differential equations for $\alpha'(x_1,x_2)$, where $\phi\in(0,\infty)$, is integrable:
\begin{gather}
\alpha'_{x_1}+\alpha_{x_2}=\frac{2}{\sinh\phi}\sinh\biggl(\frac{\alpha'}{2}\biggr)\cosh\Bigl(\frac{\alpha}{2}\Bigr)
+\frac{2\tau\cosh\phi}{\sinh\phi}\cosh\biggl(\frac{\alpha'}{2}\biggr)\sinh\Bigl(\frac{\alpha}{2}\Bigr),\nonumber\\
\alpha'_{x_2}+\alpha_{x_1}=-\frac{2}{\sinh\phi}\cosh\biggl(\frac{\alpha'}{2}\biggr)\sinh\Bigl(\frac{\alpha}{2}\Bigr)
-\frac{2\tau\cosh\phi}{\sinh\phi}\sinh\biggl(\frac{\alpha'}{2}\biggr)\cosh\Bigl(\frac{\alpha}{2}\Bigr).\label{Backcase3}
\end{gather}
Moreover, $\alpha'$ is also a solution of \eqref{eq:SHGmin}.
If $\alpha'$ and $\alpha''$ are solutions of \eqref{eq:SHGmin} obtained by the B\"acklund transformation \eqref{Backcase3}, for constants $\phi_1\neq \phi_2$, then a fourth solution $\alpha^*$ of \eqref{eq:SHGmin} is obtained algebraically with
 the superposition formula
\begin{equation*}
\tanh\biggl(\frac{\alpha^*-\alpha}{4}\biggr) =
\tau \frac{\sinh\bigl(\frac{\phi_2+\phi_1}{2}\bigr) } {\sinh\bigl(\frac{\phi_2-\phi_1}{2}\bigr)}
\tanh\biggl(\frac{\alpha'-\alpha''}{4}\biggr).
\end{equation*}
\end{Theorem}

{\bf Case 4.}
Considering $s=1$, $\epsilon=-1$, $r=1$ and $\delta=1$ then $l=-1$ and \eqref{eq:lambdaephi} implies that
$\Lambda^2 +\lambda^2=1$, hence $\Lambda=\cos\phi$ and $\lambda=\sin\phi$,
where $\phi\in (0,\pi)$.

\begin{Theorem}\label{Sinhgordon}
For any solution $\alpha(x_1,x_2)$ of
\begin{equation}\label{eq:SHG}
\alpha_{x_1x_1}-\alpha_{x_2x_2}=\sinh \alpha,
\end{equation}
 the following system of differential equations for $\alpha'(x_1,x_2)$, where $\phi\in(0,\pi)$, is integrable:
\begin{gather}
\alpha'_{x_1}+\alpha_{x_2}=\frac{2}{\sin\phi}\sinh\biggl(\frac{\alpha'}{2}\biggr)\cosh\Bigl(\frac{\alpha}{2}\Bigr)
+\frac{2\tau\cos\phi}{\sin\phi} \cosh\biggl(\frac{\alpha'}{2}\biggr)\sinh\Bigl(\frac{\alpha}{2}\Bigr),\nonumber\\
\alpha'_{x_2}+\alpha_{x_1}=-\frac{2}{\sin\phi}\cosh\biggl(\frac{\alpha'}{2}\biggr)\sinh\Bigl(\frac{\alpha}{2}\Bigr)
-\frac{2\tau\cos\phi}{\sin\phi}\sinh\biggl(\frac{\alpha'}{2}\biggr)\cosh\Bigl(\frac{\alpha}{2}\Bigr).\label{Backcase4}
\end{gather}
 Moreover, $\alpha'$ is also a solution of \eqref{eq:SHG}.
If $\alpha'$ and $\alpha''$ are solutions of \eqref{eq:SHG} obtained by the B\"acklund transformation \eqref{Backcase4}, for constants $\phi_1\neq \phi_2$, then a fourth solution $\alpha^*$ of \eqref{eq:SHG} is obtained algebraically with
 the superposition formula
\begin{equation}\label{superposcase4}
\tanh\biggl(\frac{\alpha^*-\alpha}{4}\biggr) =
\tau \frac{\sin\bigl(\frac{\phi_2+\phi_1}{2}\bigr) } {\sin\bigl(\frac{\phi_2-\phi_1}{2}\bigr)}
\tanh\biggl(\frac{\alpha'-\alpha''}{4}\biggr).
\end{equation}
\end{Theorem}

As in the first two cases, we observe that Cases 3 and 4 deal with the sinh-Gordon equation and they are complementary in the following sense: in Case 3 the parameter
 $\phi\in(0,\infty)$ and hence $\Lambda\in (1,+\infty)$, while in Case 4 we have $\phi\in (0,\pi)$ and hence
$\Lambda\in (-1,1)$.

{\bf Case 5.}
 Considering $s=1$, $\epsilon=1$, $r=0$ and $\delta=-1$ then $l=-1$ and \eqref{eq:lambdaephi} implies that
$-\Lambda^2 +\lambda^2=1$, hence $\Lambda=\sinh\phi$ and $\lambda=\cosh\phi$,
where $\phi\in [0,\infty)$.

\begin{Theorem}\label{elipticSinhgordon}
 For any solution $\alpha(x_1,x_2)$ of the elliptic sinh-Gordon equation
\begin{equation}\label{elsinhg}
\alpha_{x_1x_1}+\alpha_{x_2x_2}=\sinh \alpha,
\end{equation}
 the following system of differential equations for $\alpha'(x_1,x_2)$, where $\phi\in [0,\infty)$, is integrable:
\begin{gather}
\alpha'_{x_1}-\alpha_{x_2}=\frac{2}{\cosh\phi}\sin\biggl(\frac{\alpha'}{2}\biggr)\cosh\Bigl(\frac{\alpha}{2}\Bigr)+\frac{2\tau\sinh\phi}{\cosh\phi} \cos\biggl(\frac{\alpha'}{2}\biggr)\sinh\Bigl(\frac{\alpha}{2}\Bigr),\nonumber\\
\alpha'_{x_2}+\alpha_{x_1}=-\frac{2}{\cosh\phi}\cos\biggl(\frac{\alpha'}{2}\biggr)\sinh\Bigl(\frac{\alpha}{2}\Bigr)+\frac{2\tau\sinh\phi}{\cosh\phi} \sin\biggl(\frac{\alpha'}{2}\biggr)\cosh\Bigl(\frac{\alpha}{2}\Bigr).\label{BTcase5}
\end{gather}
 Moreover, $\alpha'$ is a solution of the elliptic sine-Gordon equation
\begin{equation}\label{elipticSG}
\alpha'_{x_1x_1}+\alpha'_{x_2x_2}=\sin \alpha'.
\end{equation}
If $\alpha'$ and $\alpha''$ are solutions of \eqref{elipticSG} obtained by the B\"acklund transformation \eqref{BTcase5}, for constants $\phi_1\neq \phi_2$, then a solution $\alpha^*$ of \eqref{elsinhg} is obtained algebraically with the superposition formula
\begin{equation}\label{superposcase5}
\tanh \frac{\alpha^*}{2}=
 \frac{P-Q\tanh {\frac{\alpha}{2}}} {Q-P \tanh {\frac{\alpha}{2}}},
\end{equation}
where
\[ P=\tau(\sinh\phi_2-\sinh\phi_1) \sin\biggl(\frac{\alpha''-\alpha'}{2}\biggr),
\]
\[
Q=(1+\sinh\phi_1\sinh\phi_2)\cos\biggl(\frac{\alpha''-\alpha'}{2}\biggr) -\cosh\phi_1\cosh\phi_2.
\]
\end{Theorem}

{\bf Case 6.}
Finally, considering $s=1$, $\epsilon=1$, $r=1$ and $\delta=-1$ then $l=1$ and \eqref{eq:lambdaephi} implies that
$-\Lambda^2 +\lambda^2=1$, hence $\Lambda=\sinh\phi$ and $\lambda=\cosh\phi$,
where $\phi\in [0,\infty)$.

\begin{Theorem}\label{elipticSinegordon}
 For any solution $\alpha(x_1,x_2)$ of the elliptic sine-Gordon equation \eqref{elipticSG},
the following system of differential equations for $\alpha'(x_1,x_2)$, where $\phi\in [0,\infty)$, is integrable:
\begin{gather}
\alpha'_{x_1}-\alpha_{x_2}=\frac{2}{\cosh\phi}\sinh\biggl(\frac{\alpha'}{2}\biggr)\cos\Bigl(\frac{\alpha}{2}\Bigr)+\frac{2\tau\sinh\phi}{\cosh\phi} \cosh\biggl(\frac{\alpha'}{2}\biggr)\sin\Bigl(\frac{\alpha}{2}\Bigr),\nonumber\\
\alpha'_{x_2}+\alpha_{x_1}=-\frac{2}{\cosh\phi}\cosh\biggl(\frac{\alpha'}{2}\biggr)\sin\Bigl(\frac{\alpha}{2}\Bigr)+\frac{2\tau\sinh\phi}{\cosh\phi} \sinh\biggl(\frac{\alpha'}{2}\biggr)\cos\Bigl(\frac{\alpha}{2}\Bigr).\label{BTcase6}
\end{gather}
 Moreover, $\alpha'$ is a solution of the elliptic sinh-Gordon equation \eqref{elsinhg}. If $\alpha'$ and $\alpha''$ are solutions of \eqref{elsinhg} obtained by the B\"acklund transformation \eqref{BTcase6} for constants $\phi_1\neq \phi_2$, then a solution~$\alpha^*$ of the elliptic sine-Gordon equation \eqref{elipticSG} is obtained algebraically with the superposition formula
\begin{equation}\label{superposcase6}
\tan \frac{\alpha^*}{2}= \frac{\tilde{P}-\tilde{Q}{\tan \frac{\alpha}{2} } }{\tilde{Q}+\tilde{P} {\tan\frac{\alpha}{2}} },
\end{equation}
where
\begin{gather*}
\tilde{P}=\tau(\sinh \phi_2-\sinh \phi_1) \sinh\biggl(\frac{\alpha''-\alpha'}{2}\biggr),
\\
\tilde{Q}=(1+\sinh\phi_1\sinh\phi_2)\cosh\frac{\alpha''-\alpha'}{2}-\cosh\phi_1\cosh\phi_2.
\end{gather*}
\end{Theorem}

\section{Examples}\label{examples}
In this section, we present some examples and illustrations of the previous results. Since explicit solutions of the sine-Gordon equation, the corresponding surfaces and applications of Theorem~\ref{SG}, namely the B\"acklund transformation \eqref{Backcase1} and the superposition formula \eqref{superposcase1} are well known (see, for example,~\cite{Tenenblat1998}), we will concentrate in applying the new results. Our first example provides the most straightforward application of Theorem~\ref{thBT}.

\begin{Example}\label{example1} Let $\delta=\pm 1$, $\epsilon=\pm 1$,
$r,s\in \{0,1\}$, $r\leq s$, be fixed constants satisfying \eqref{eq:rs} and \eqref{eq:epphi}. Consider $\alpha=0$ a trivial solution of \eqref{eqorig}. Then the B\"acklund transformation \eqref{eq:btanalytic_v} reduces to
\begin{gather*}
\alpha'_{x_1}=\frac{2}{\lambda}S_{(-1)^r} \biggl(\frac{\alpha'}{2}\biggr),\qquad
\alpha'_{x_2}=-\frac{2\delta\tau\Lambda}{\lambda} S_{(-1)^r}\biggl(\frac{\alpha'}{2}\biggr),
\end{gather*}
where $\tau=\pm 1$, and the constants $\Lambda$, $\lambda$ satisfy \eqref{eq:lambdaephi}. The solutions of this system are given by
\begin{equation}\label{eq:BToftrivialsolution}
\alpha'(x_1,x_2)=\begin{cases}
4\arctan \exp\xi(x_1,x_2)&\mbox{if} \ r=0,\\
2\ln\tanh\left(-\frac{1}{2}\xi(x_1,x_2)\right)&\mbox{if} \ r=1,
\end{cases}
\end{equation}
where
\[
\xi(x_1,x_2)=\frac{1}{\lambda}x_1-\delta\tau\frac{\Lambda}{\lambda}x_2+c,
\]
and $c$ is a constant that depends on a given initial condition.
Moreover, $\alpha'$ is a solution of \eqref{eqprime}.
\end{Example}

Although $\alpha=0$ does not correspond to any surface, we observe that $\alpha'$ does correspond to a surface (see, for instance, \eqref{superfalpha1caso4} in Example~\ref{example2}).

In the following examples, we include explicit solutions of the differential equations
and some of the surfaces of constant curvature that correspond to those solutions.
We recall from~\cite{KRT2022}, that under the hypothesis of Theorem~\ref{th:integrabilitytheorem}, if $X(x_1,x_2)$ is a parametrized surface in ${\mathbb{R}}^3_s$ of curvature $\delta=\pm 1$, then there exists a two-parameter family of surfaces $\bar{X}$
related to $X$ by a B\"acklund-type line congruence. Moreover, it follows from Theorem~\ref{thBT} that if
$\alpha$ is a solution of \eqref{eqorig}, that corresponds to the surface $X$, then the B\"acklund transformation \eqref{eq:btanalytic_v} is integrable for $\alpha'$ and the surface $\bar{X}$ that corresponds to $\alpha'$ is given by
\begin{equation*}
\bar{X}=X+\lambda
\Biggl(
\frac{C_{(-1)^r}\bigl(\frac{\alpha'}{2}\bigr)}{C_l\bigl(\frac{\alpha}{2}\bigr)}X_{x_1}
+
\frac{S_{(-1)^r}\bigl(\frac{\alpha'}{2}\bigr)}{S_l\bigl(\frac{\alpha}{2}\bigr)}X_{x_2}
\Biggr).
\end{equation*}

\begin{figure}
 \centering
 \begin{minipage}{.31\linewidth}
 \centering
 {\includegraphics[height=110pt]{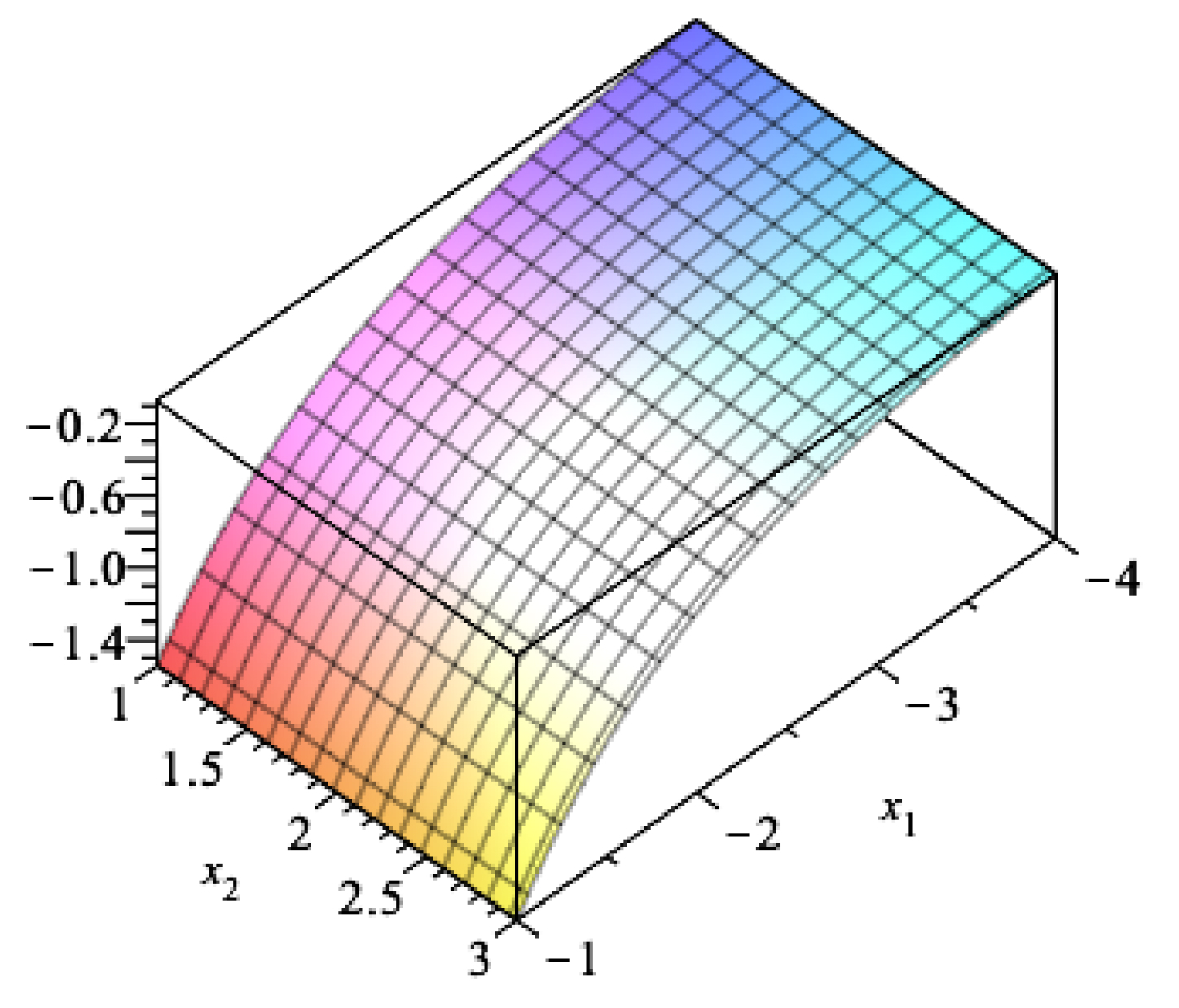}}\\
 $\alpha_1$
 \end{minipage}
 \begin{minipage}{.31\linewidth}
 \centering
 {\includegraphics[height=110pt]{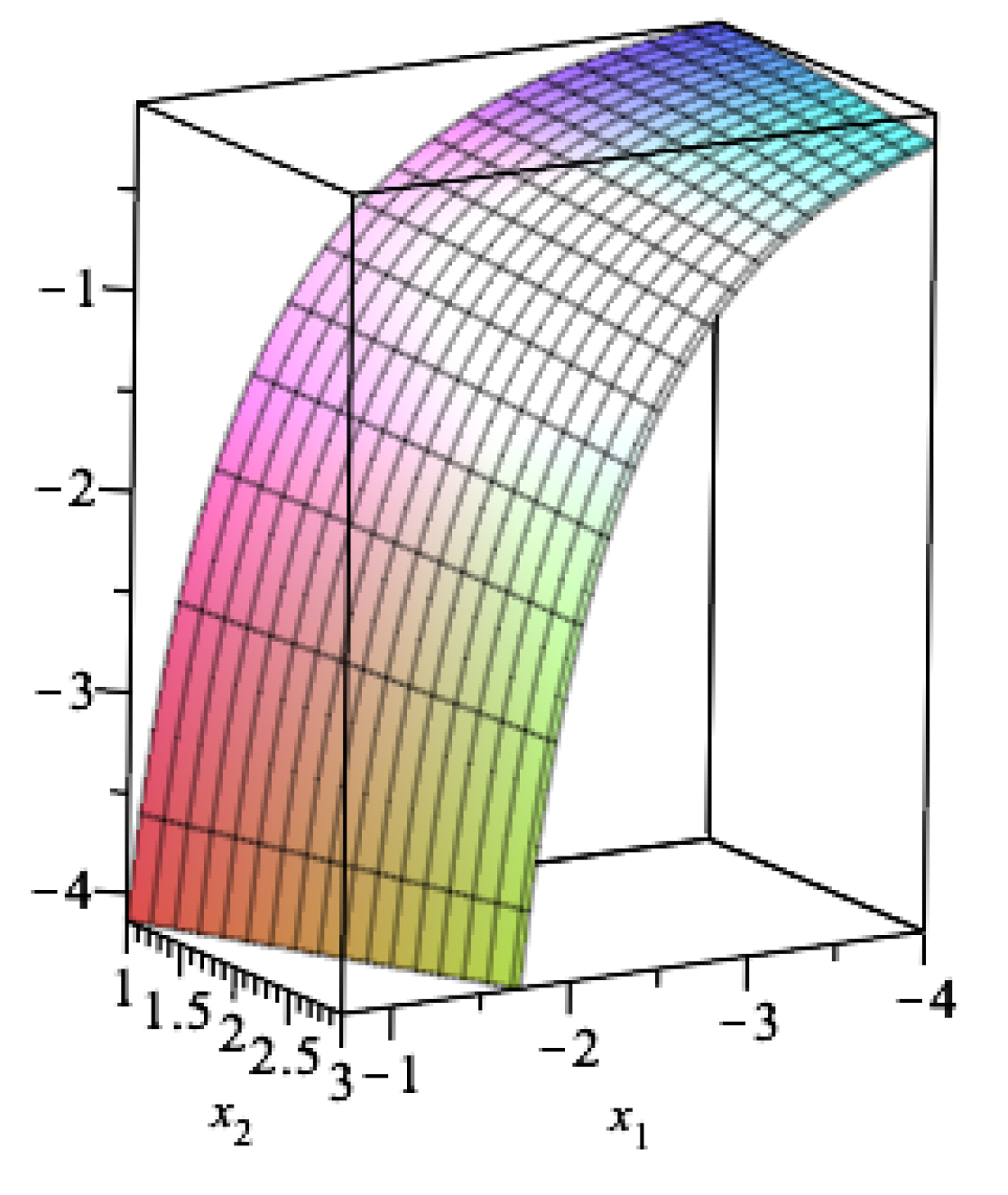}}\\
 $\alpha_2$
 \end{minipage}
 \begin{minipage}{.31\linewidth}
 \centering
 {\includegraphics[height=110pt]{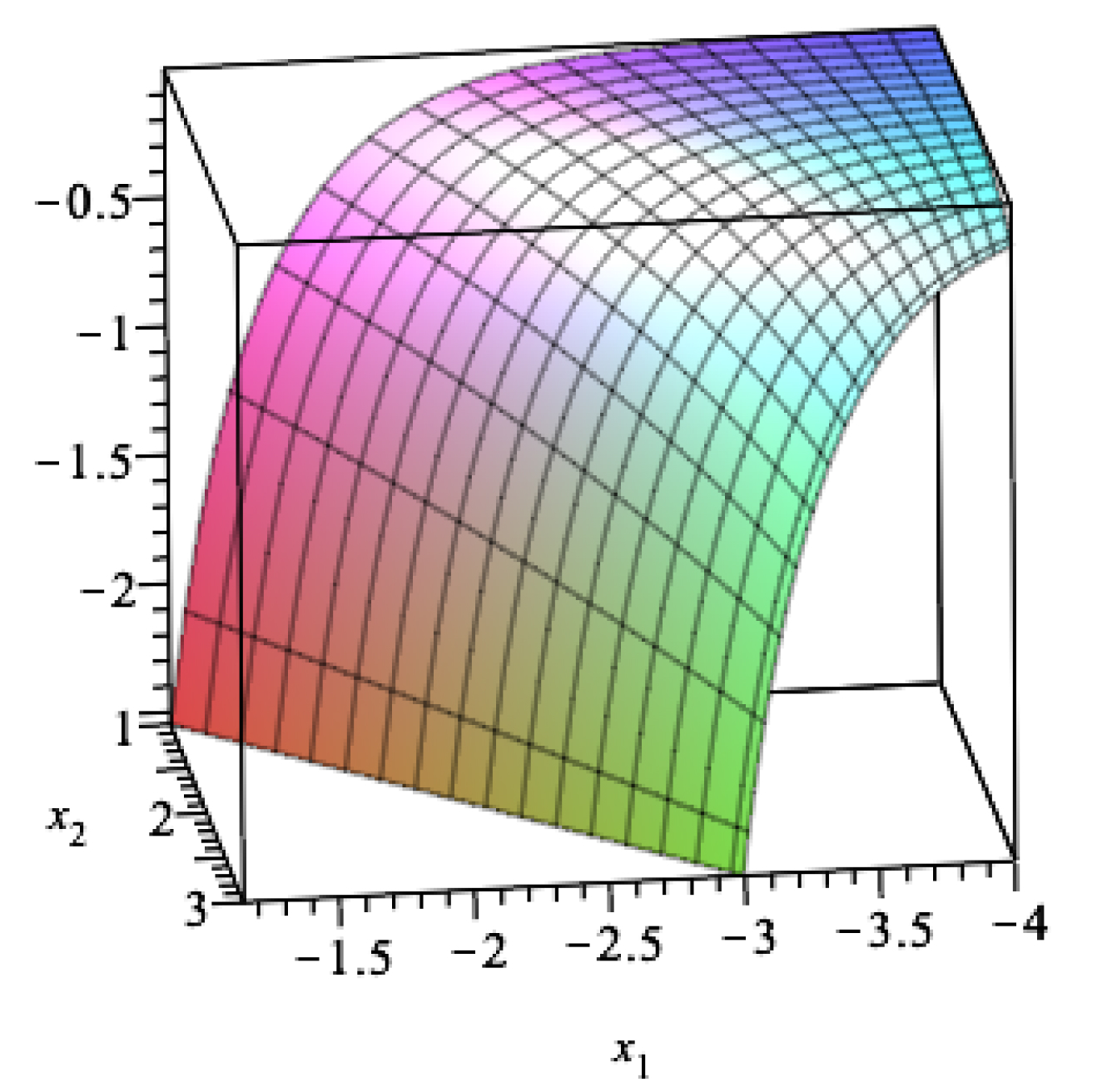}}\\
 $\alpha_3$
 \end{minipage}

 \begin{minipage}{.31\linewidth}
 \centering
 {\includegraphics[height=110pt]{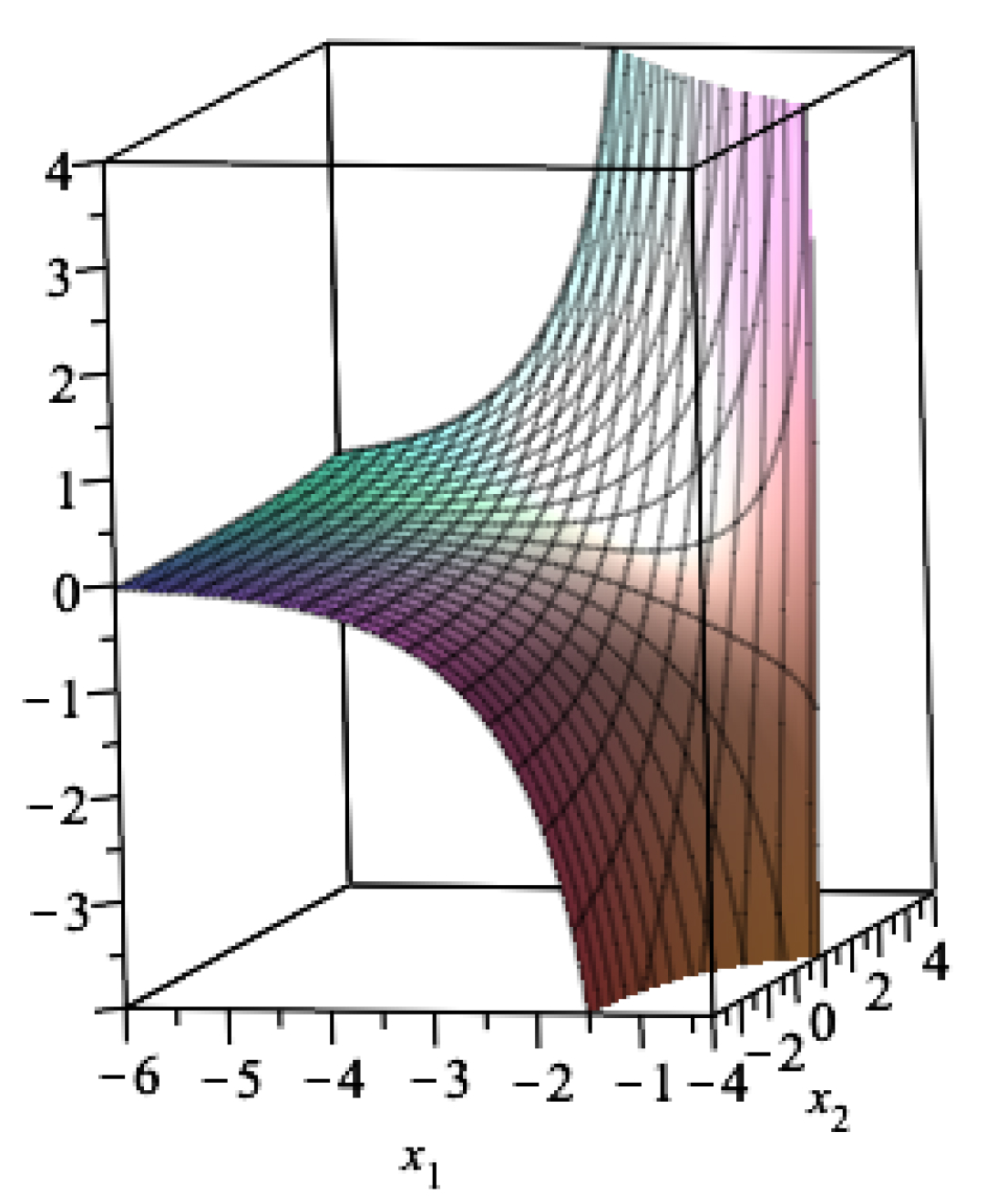}}\\
 $\alpha_{12}$
 \end{minipage}
 \begin{minipage}{.31\linewidth}
 \centering
 {\includegraphics[height=110pt]{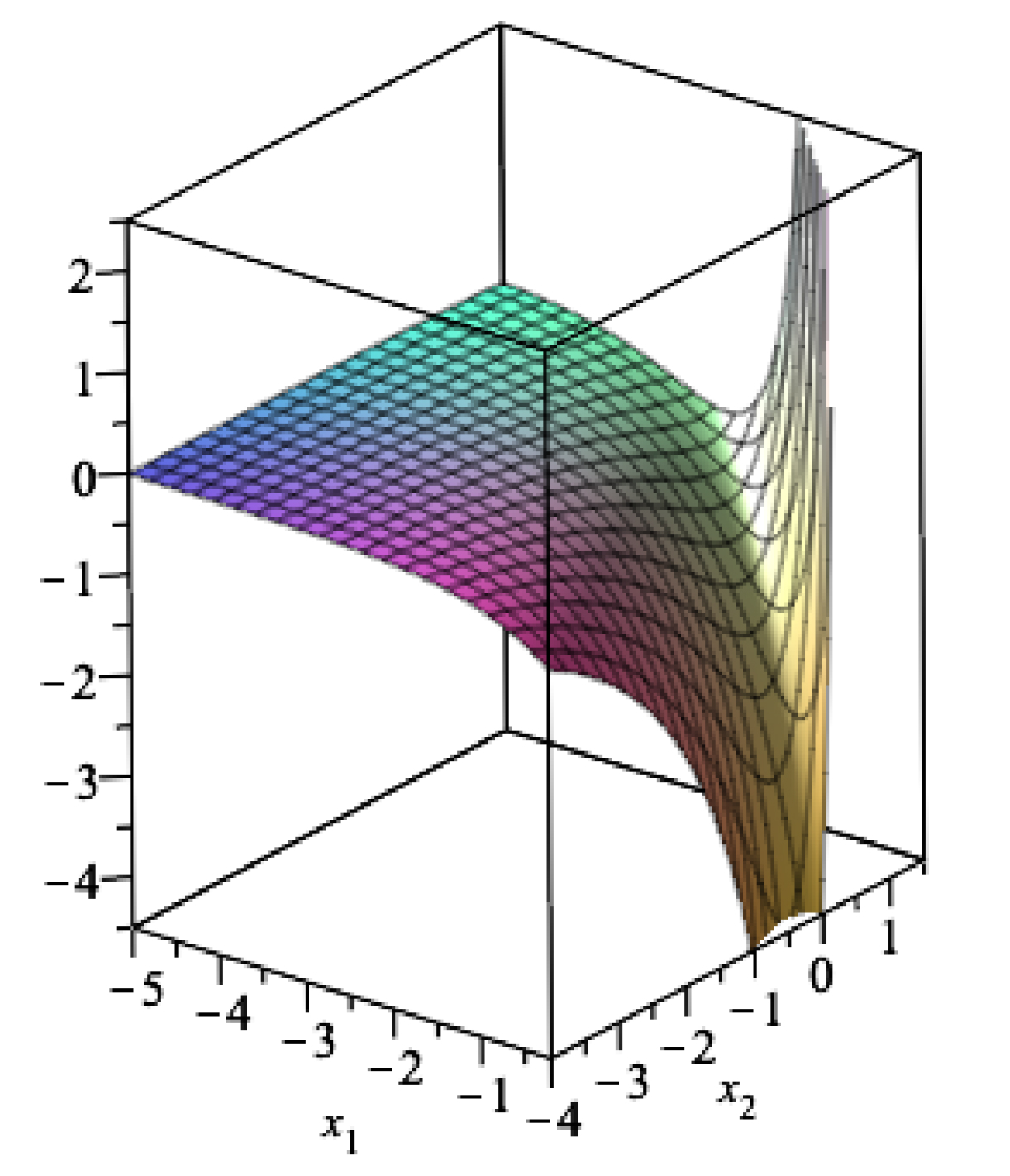}} \\
 $\alpha_{23}$
 \end{minipage}
 \begin{minipage}{.31\linewidth}
 \centering
 {\includegraphics[height=110pt]{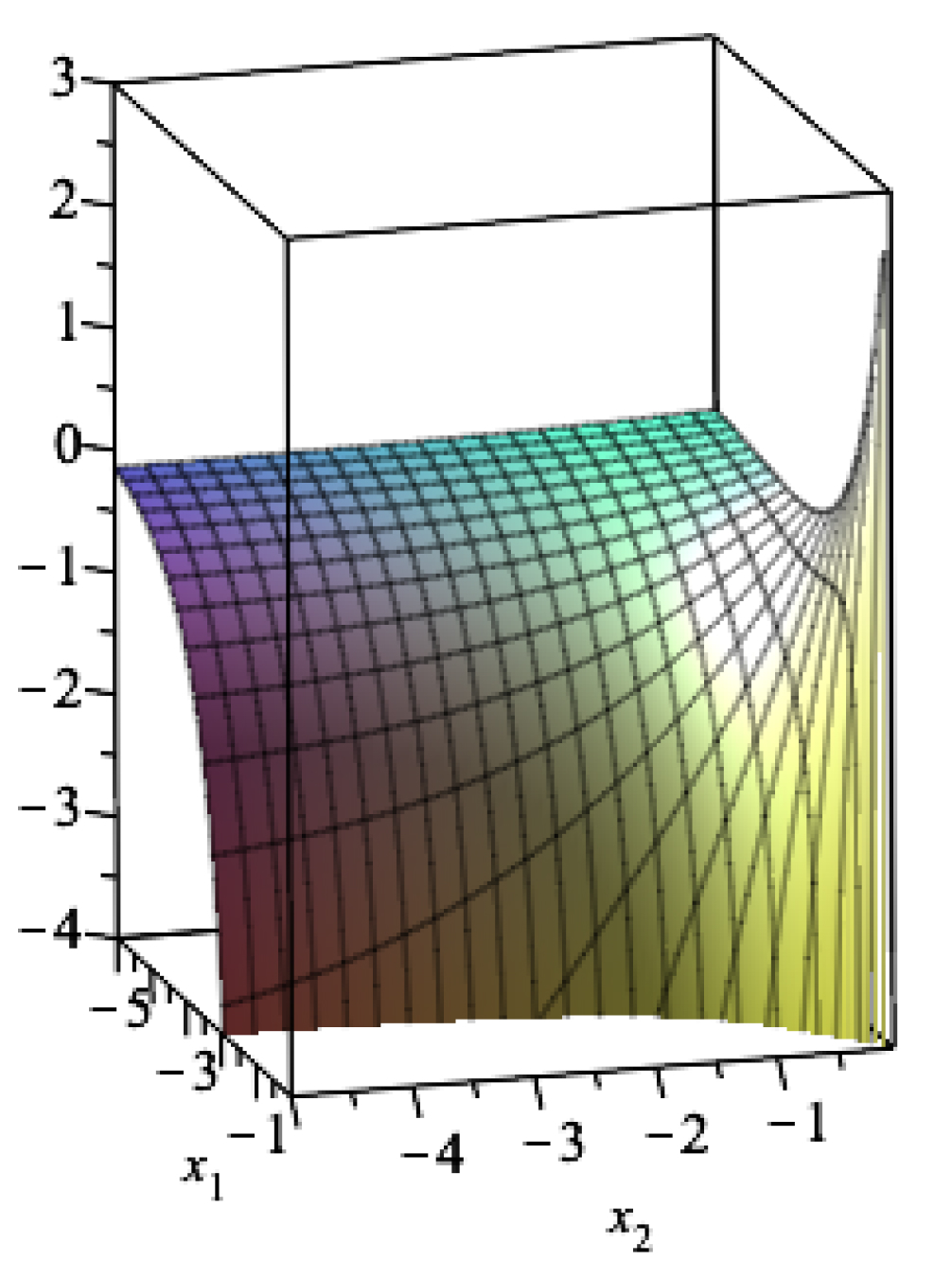}}\\
 $\alpha_{123}$
 \end{minipage}
 \caption {Solutions of the sinh-Gordon equation: $\alpha_1$, $\alpha_2$ and $\alpha_3$ are obtained from the null solution by the B\"acklund transformation with parameters $\phi_1=\pi/2$, $\phi_2=\pi/3$ and $\phi_3=\pi/8$ respectively. $\alpha_{12}$~is a~solution obtained from $\alpha_1$ and $\alpha_2$ by the superposition formula. $\alpha_{23}$ is a~solution obtained from $\alpha_2$ and~$\alpha_3$ by the superposition formula and $\alpha_{123}$ is obtained from $\alpha_{12}$ and~$\alpha_{23}$ by the superposition formula. \label{fig:caso4alpha}}
 \end{figure}

\begin{figure}
 \centering
 \begin{minipage}{.32\linewidth}
 \centering
 {\includegraphics[height=130pt]{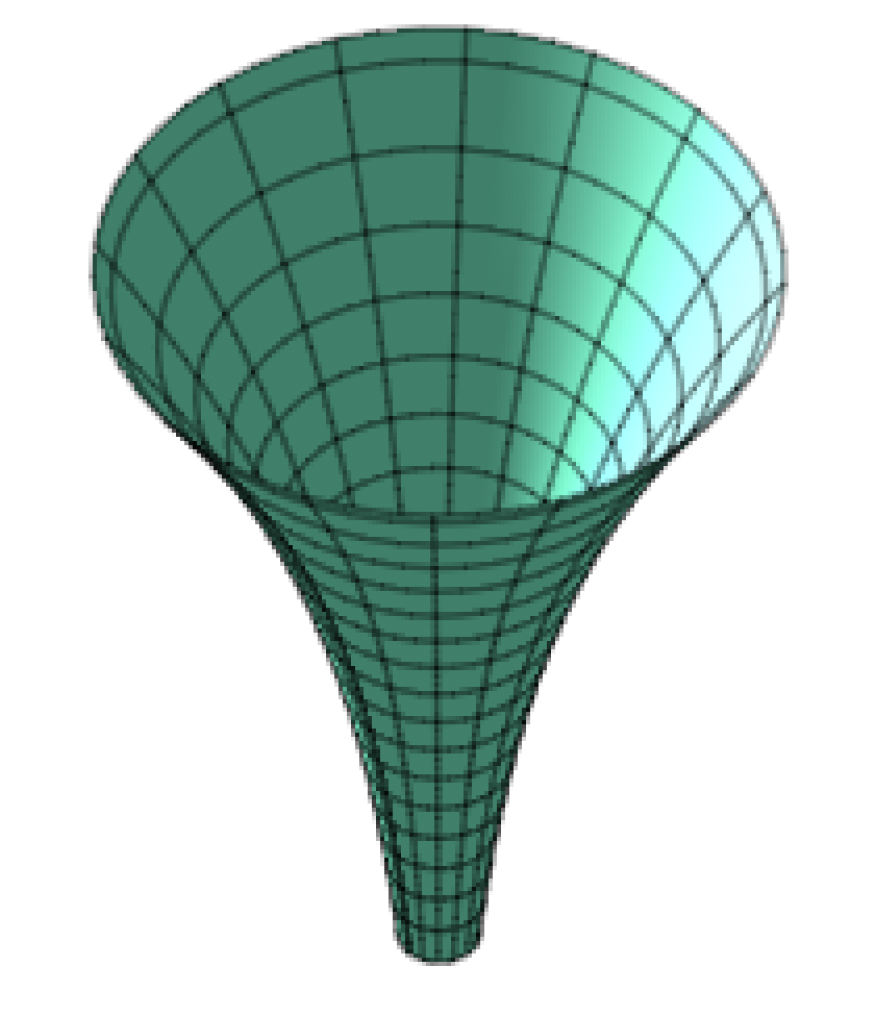}}\\
 $X_1$
 \end{minipage}
 \begin{minipage}{.32\linewidth}
 \centering
 {\includegraphics[height=130pt]{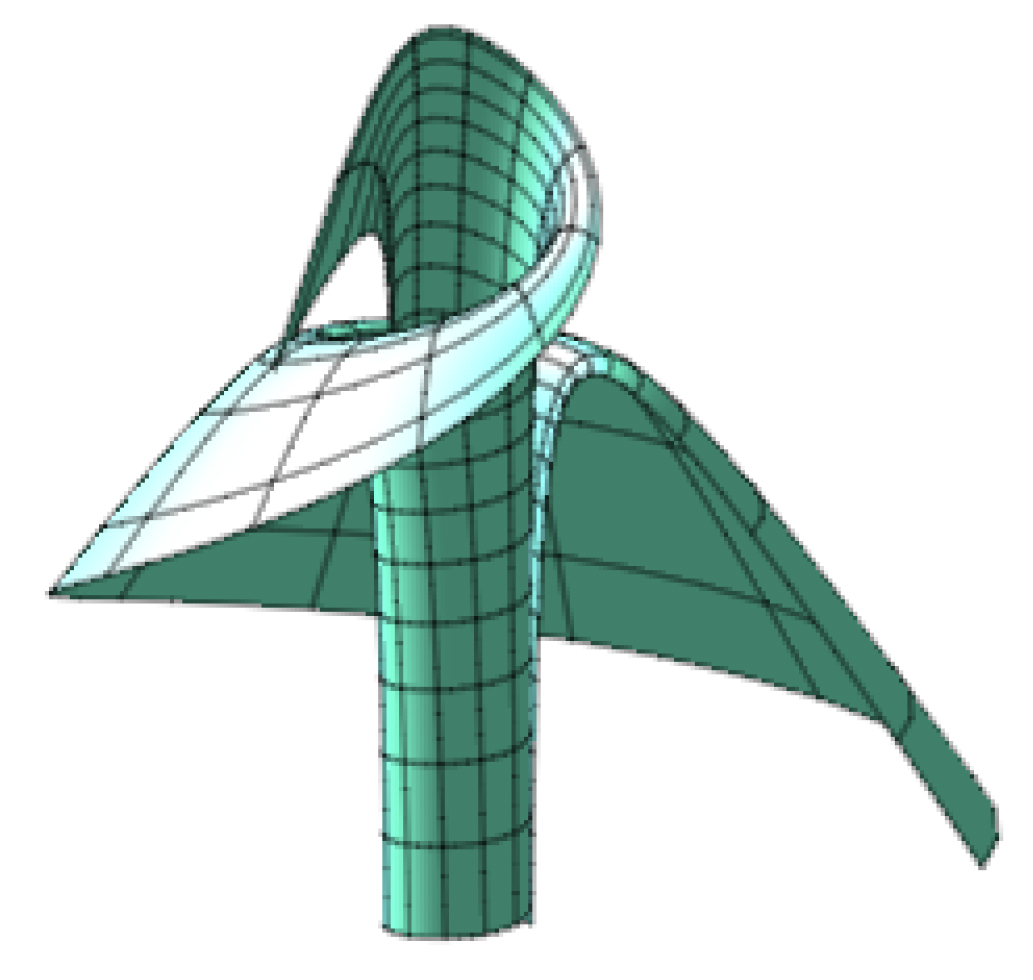}}\\ 
 $X_{12}$
\end{minipage}
 \begin{minipage}{.32\linewidth}
 \centering
 {\includegraphics[height=130pt]{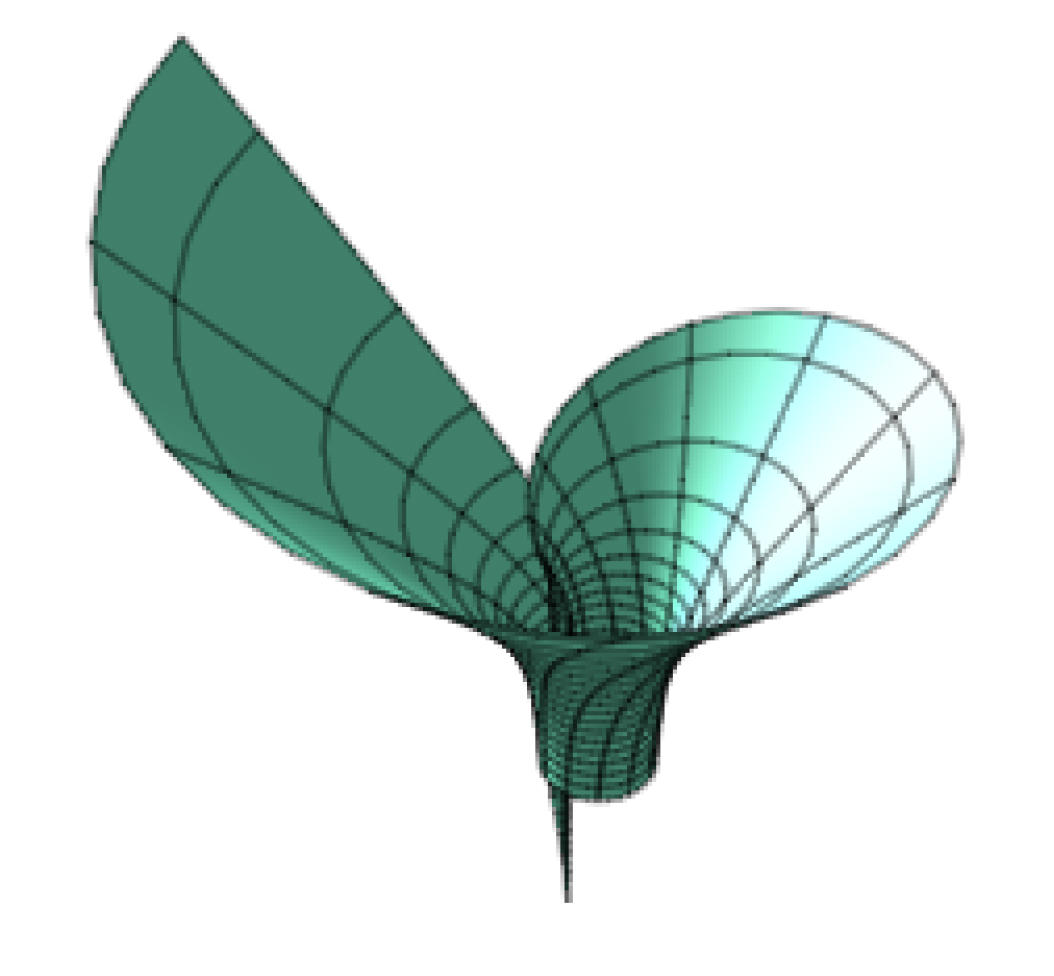}}\\
 $X_{123}$
 \end{minipage}
 \caption {$X_1$ is a time-like surface with Gaussian curvature $K=1$ contained in ${\mathbb{R}}^3_1$ corresponding to the solution $\alpha_1$. $X_{12}$ is a time-like surface obtained from the surface $X_1$ by a B\"acklund transformation, with parameter $\phi_2=\pi/3$ and it corresponds to the solution $\alpha_{12}$.
 $X_{123}$ is the surface obtained from $X_{12}$ by a B\"acklund transformation, with parameter $\phi_3=\pi/8$ and it corresponds to $\alpha_{123}$.
 \label{fig:caso4X}}
\end{figure}

\begin{Example}\label{example2}
 We consider the constants $\delta=1$, $\epsilon=-1$, $r=s=1$,
then \eqref{eqorig} is the sinh-Gordon equation (SHGE). We will illustrate Theorem \ref{Sinhgordon}, by starting with the trivial solution $\alpha=0$. From Example \ref{example1}, the B\"acklund transformation \eqref{Backcase4} provides new solutions for the SHGE given by \eqref{eq:BToftrivialsolution}, with $r=1$.
 Let $\alpha_j=2\ln \tanh(-\xi_j(x_1,x_2)/2)$, $j=1,2,3$, be a solution of the SHGE given by
$\xi_j=1/\sin(\phi_j)(x_1+\cos(\phi_j)x_2)+c_j$, for distinct parameters $\phi_j\in(0,\pi)$, $\tau=-1$ and arbitrary constants $c_j$.
Observe that for each $j$, $\alpha_j$ is defined on the half-plane $\xi_j<0$. We will apply the superposition formula to distinct pairs of $\alpha_j$.

Now we choose $\phi_1=\pi/2$, $\phi_2=\pi/3$ and $\phi_3=\pi/8$. By applying
the superposition formula~\eqref{superposcase4} to $\alpha$ and $\alpha_j$, $j=1,2$ and the parameters $\phi_1$, $\phi_2$ (here we are denoting the solutions $\alpha'$, $\alpha''$ and $\alpha^*$ by $\alpha_1$, $\alpha_2$ and $\alpha_{12}$ respectively), one gets the solution
$\alpha_{12}$.
Similarly, by applying the superposition formula to $\alpha$, $\alpha_2$ and $\alpha_3$ with the parameters $\phi_2$ and $\phi_3$, one gets the solution $\alpha_{23}$. These solutions of the sinh-Gordon equation are given by
\begin{gather*}
\alpha_{12}=4\operatorname{arctanh}\Biggl(\frac{\sin\bigl(\frac{5\pi}{12}\bigr)}{\sin\bigl(\frac{\pi}{12}\bigr)}
\tanh\biggl(\frac{\alpha_1-\alpha_2}{4}\biggr)\Biggr),
\\
\alpha_{23}=4\operatorname{arctanh}\Biggl(\frac{\sin\bigl(\frac{11\pi}{48}\bigr)}{\sin\bigl(\frac{5\pi}{48}\bigr)}
\tanh\biggl(\frac{\alpha_2-\alpha_3}{4}\biggr)\Biggr),
\end{gather*}
 where $\tanh((\alpha_1-\alpha_2)/{4})$ and $\tanh((\alpha_2-\alpha_3)/{4})$ are functions of $(x_1,x_2)$ determined by equation $\exp(\alpha_j/2)=\tanh(-\xi_j/2)$. Moreover, the solution $\alpha_{12}$ tends to $\pm \infty$ on the two curves given by $\tanh(\alpha_1-\alpha_2)/4=\pm \bigl(\sin\bigl(\frac{\pi}{12}\bigr)/\sin\bigl(\frac{5\pi}{12}\bigr)\bigr)$. Similarly, $\alpha_{23}$ tends to $\pm \infty$ on two curves. The superposition of~$\alpha_{12}$ and $\alpha_{23}$ with the parameters $\phi_1$ and $\phi_3$, provides the solution $\alpha_{123}$ (see diagram~\eqref{diagram}),
\[
\alpha_{123}= \alpha_2 + 4\operatorname{arctanh}\Biggl(\frac{\sin\bigl(\frac{5\pi}{16}\bigr)}{\sin\bigl(\frac{3\pi}{16}\bigr)}
\tanh\biggl(\frac{\alpha_{12}-\alpha_{23}}{4}\biggr)\Biggr).
\]
The graph of the functions $\alpha_1$, $\alpha_2$, $\alpha_3$ and $\alpha_{12}$, $\alpha_{23}$, $\alpha_{123}$, with the constants $c_j=0$, can be visualized in Figure~\ref{fig:caso4alpha}.
This procedure provides infinitely many solutions of the sinh-Gordon equation, whose explicit expressions may get more and more complicated.

The time-like surface $X_1$ in ${\mathbb{R}}^3_1$, with constant curvature $K=1$, that corresponds to the function $\alpha_1$, is given by
\begin{equation}\label{superfalpha1caso4}
X_1=\biggl( \frac{\cos x_2}{\sinh x_1}, \frac{\sin x_2}{\sinh x_1}, x_1-\frac{\cosh x_1}{\sinh x_1} \biggr).
\end{equation}
 The surface $X_{12}$ corresponding to $\alpha_{12}$ is given
by
\[
X_{12}= X_1+\sin(\phi_2)\Biggl( \frac{ \cosh\bigl(\frac{\alpha_{12}}{2}\bigr)}{\sinh\bigl(\frac{\alpha_1}{2}\bigr)} X_{1,x_1} + \frac{ \sinh\bigl(\frac{\alpha_{12}}{2}\bigr)}{\sinh\bigl(\frac{\alpha_1}{2}\bigr)} X_{1,x_2}\Biggr).
\]
The solution $\alpha_{123}$ corresponds to a surface $X_{123}$,
\[
X_{123}= X_{12}+\sin\phi_3\Biggl( \frac{ \cosh\bigl(\frac{\alpha_{123}}{2}\bigr)}{\sinh\bigl(\frac{\alpha_{12}}{2}\bigr)} X_{12,x_1} + \frac{ \sinh\bigl(\frac{\alpha_{123}}{2}\bigr)}{\sinh\bigl(\frac{\alpha_{12}}{2}\bigr)} X_{12,x_2}\Biggr).
\]
These surfaces are illustrated in Figure~\ref{fig:caso4X}. We point out that although the surface $X_1$ in Figure~\ref{fig:caso4X} looks like a pseudosphere, considering this surface with the metric induced by the Euclidean space ${\mathbb{R}}^3$, its Gaussian curvature is given by $-\sinh^4(x_1)/ \bigl(1+\cosh^2(x_1)\bigr)^2$, which is not constant. Moreover, we observe that the surfaces in this example were obtained due to the superposition formulae.
\end{Example}


\begin{Example}
 We consider the constants $\delta=-1$, $\epsilon=1$, $r=0$, $s=1$,
then \eqref{eqorig} is the elliptic sinh-Gordon equation (ESHGE). We will illustrate Theorem \ref{elipticSinhgordon}, by starting with the trivial solution $\alpha=0$. The B\"acklund transformation \eqref{BTcase5} provides solutions for the elliptic sine-Gordon equation (ESGE) given by \eqref{eq:BToftrivialsolution}, where $r=0$,
 $\Lambda=\sinh\phi$, $\lambda=\cosh\phi$ and any \mbox{$\phi\in[0,+\infty)$}.
Therefore, by choosing $\tau=1$ and $\phi_1\neq \phi_2$ in
\eqref{eq:BToftrivialsolution}, we have the solutions~$\alpha_1$ and~$\alpha_2$ for the ESGE given by
\begin{gather*}
\alpha_j=4 \arctan {\rm e}^{\xi_j }, \qquad
\end{gather*}
where
\begin{gather*}
\xi_j=\frac{1}{\cosh \phi_j}(x_1+\sinh\phi_jx_2)+c_j, \qquad j=1,2, \quad c_j\in{\mathbb{R}}.
\end{gather*}

 By applying the superposition formula \eqref{superposcase5}, we obtain a new solution $\alpha_{12}$ for the ESHGE,
\begin{equation}\label{alpha12caso5}
 \tanh\Bigl(\frac{\alpha_{12}}{2}\Bigr)=
 \frac{(\sinh \phi_2-\sinh \phi_1)\sin\bigl(\frac{\alpha_2-\alpha_1}{2 } \bigr) }{(1+\sinh\phi_1 \sinh\phi_2)
 \cos\bigl(\frac{\alpha_2-\alpha_1}{2 } \bigr)-\cosh\phi_1\cosh\phi_2}.
 \end{equation}
We observe that the function $\alpha_{12}$ tends to $\pm\infty$ on two curves (when the right-hand side of \eqref{alpha12caso5} tends to $\pm 1$) defined by
\begin{equation}\label{singcaso5}
\frac{{\rm e}^{\xi_2}-{\rm e}^{\xi_1 }} {1+{\rm e}^{\xi_2+\xi_1}}\pm \frac {\sinh\phi_2-\sinh\phi_1}{1+\cosh(\phi_2+\phi_1)}=0.
\end{equation}
In particular, by choosing $\phi_1=0$, $\phi_2=\ln\bigl[\bigl(1+\sqrt{5}\bigr)/2\bigr]$ and the constants $c_j=0$,
we have $\xi_1=x_1$ and $\xi_2=(2x_1+x_2)/\sqrt{5} $.
Hence, starting with the null solution of the ESHGE, the superposition formula provides
a new explicit solution of the same equation, namely
\begin{equation*}
\alpha_{12}= 2 \operatorname{arctanh}\Biggl( \frac{ \tan\bigl(\frac{\alpha_2-\alpha_1} {4 } \bigr) }
{\bigl(1-\frac{\sqrt{5}}{2}\bigr)-\bigl(1+\frac{\sqrt{5}}{2}\bigr) \bigl(\tan\bigl(\frac{\alpha_2-\alpha_1} {4 } \bigr) \bigr)^2}
 \Biggr),
\end{equation*}
where
\[
\tan\Bigl(\frac{\alpha_2-\alpha_1} {4 } \Bigr)= \frac{{\rm e}^{(2x_1+x_2)/\sqrt{5}}-{\rm e}^{x_1} } {1+ {\rm e}^{(2x_1+x_2)/\sqrt{5}+x_1} }.
\]
In Figure \ref{fig:caso5alpha}, one can visualize the graph of the solution of the ESHGE given by $\alpha_{12}$ and also the two curves \eqref{singcaso5} where it tends to $\pm\infty$.

\begin{figure}
 \centering
 \begin{minipage}{.32\linewidth}
 \centering
 {\includegraphics[height=127pt]{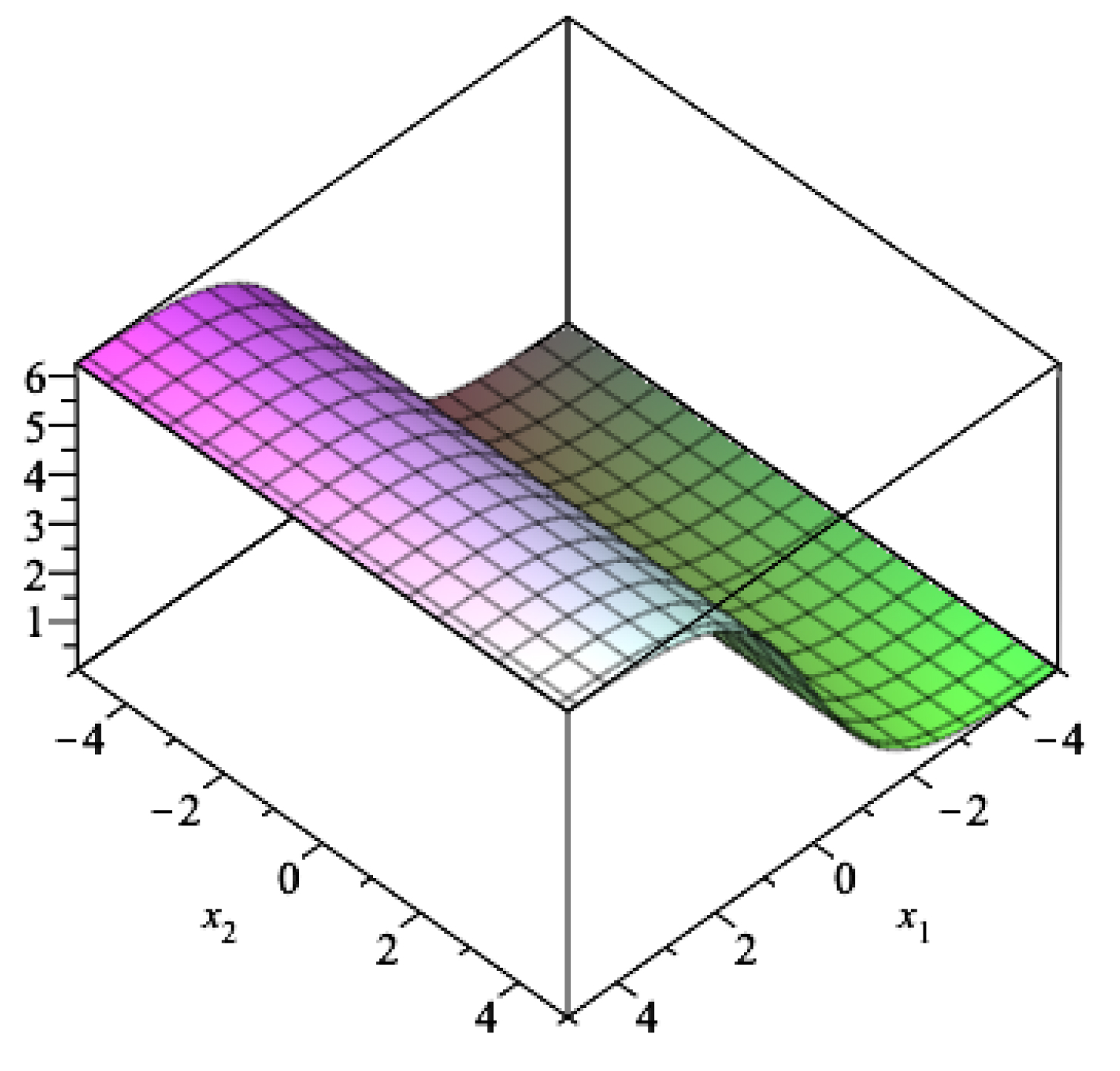}}\\
 $\alpha_1$
 \end{minipage}
 \begin{minipage}{.32\linewidth}
 \centering
 {\includegraphics[height=127pt]{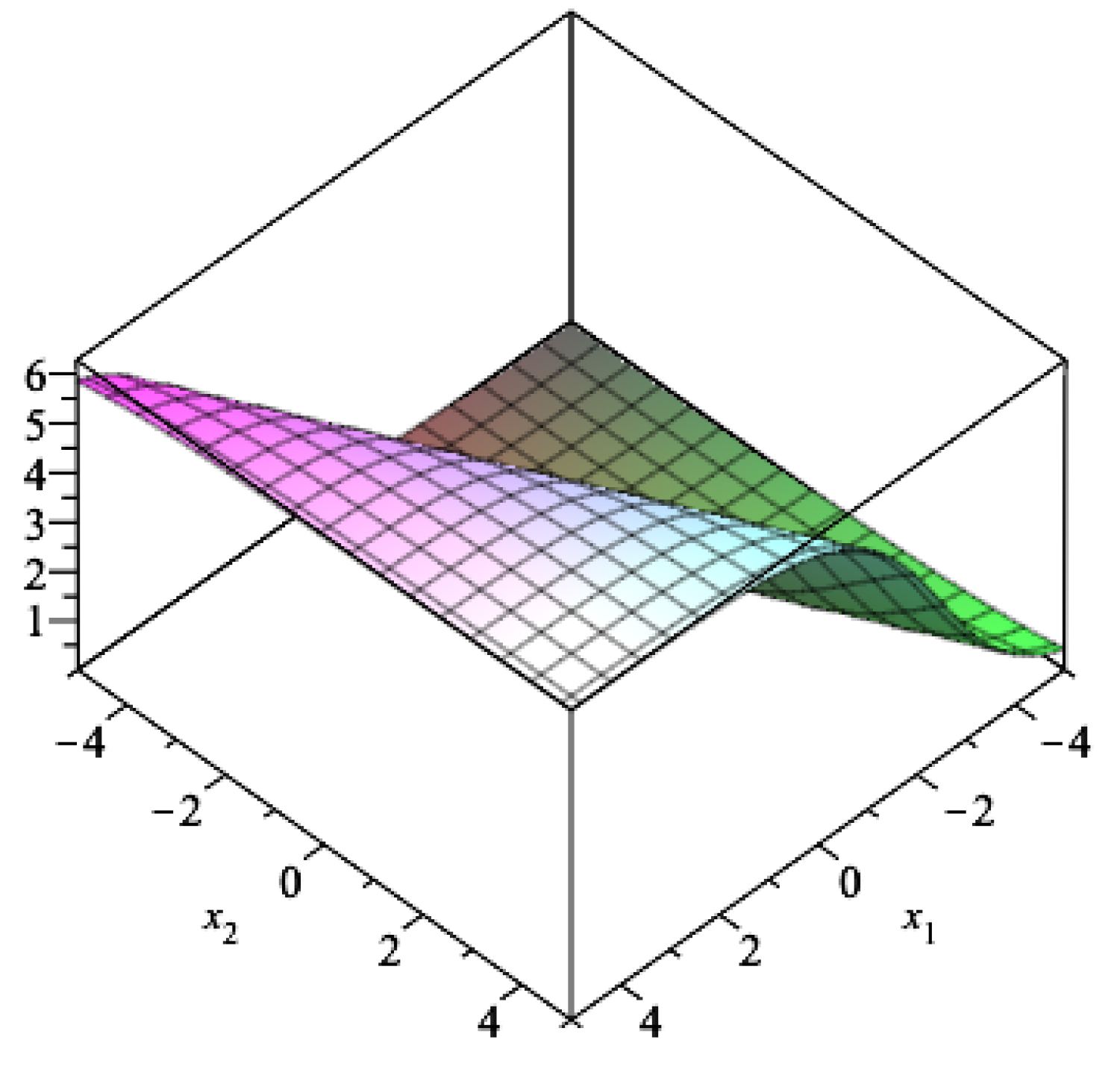}} \\
 $\alpha_2$
 \end{minipage}
 \begin{minipage}{.32\linewidth}
 \centering
 {\includegraphics[height=127pt]{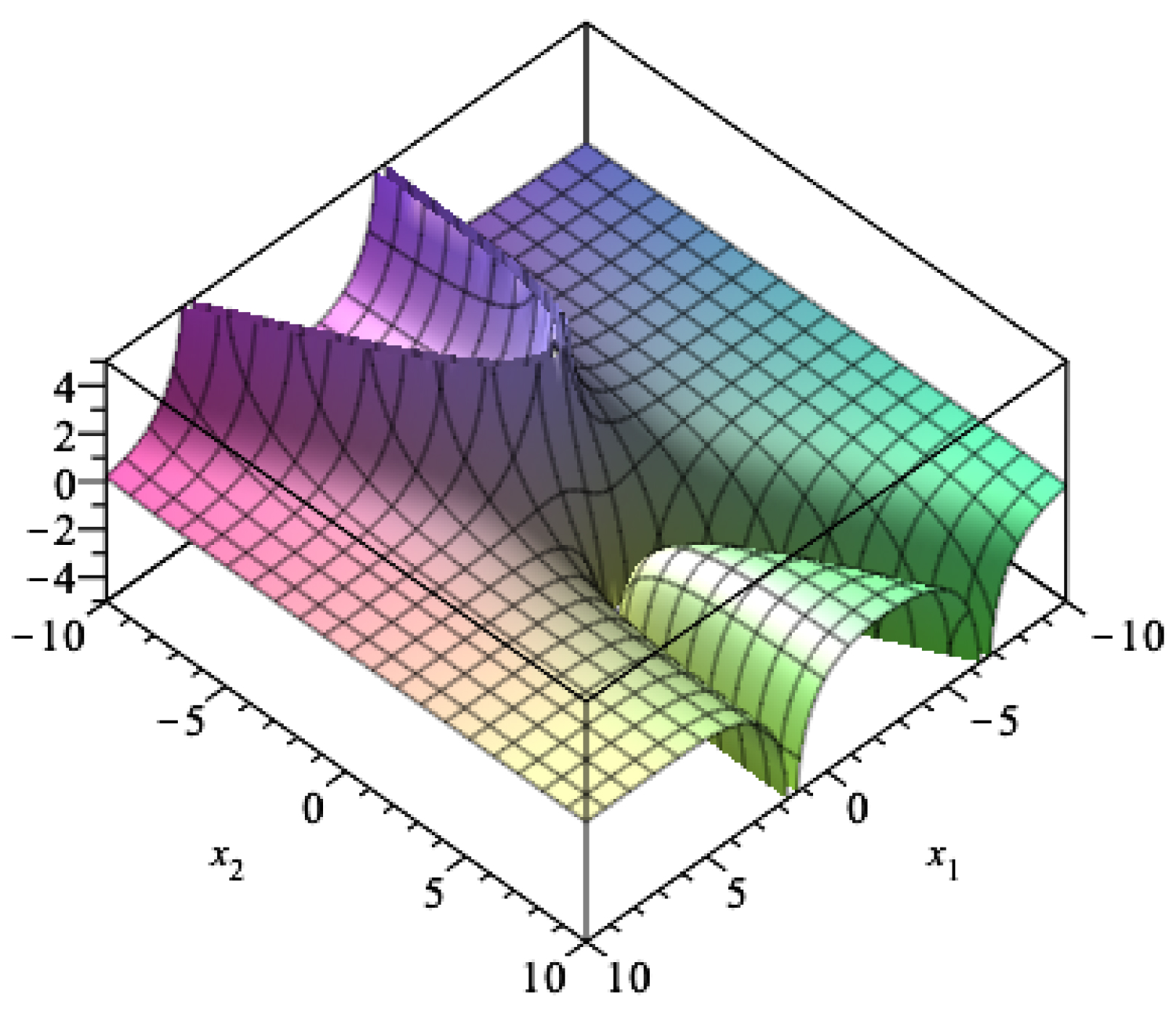}}\\
 $\alpha_{12}$
 \end{minipage}
 \caption { $\alpha_1$ and $\alpha_2$ are solutions of the elliptic sine-Gordon equation, obtained from the null solution of the elliptic sinh-Gordon equation, by B\"acklund transformations with parameters $\phi_1=0$ and $\phi_2=\ln( (1+\sqrt{5})/2)$ respectively.
 $\alpha_{12}$ is a solution of the elliptic sinh-Gordon equation defined on the complement of two curves, obtained by using the superposition formula. \label{fig:caso5alpha}}
 \end{figure}
 We observe that the surface associated to $\alpha_1$ (the intermediate step which is a solution of the ESGE) is given by \eqref{superfESGE} in the next example.
\end{Example}

\begin{Example}
In this example, we illustrate Theorem \ref{elipticSinegordon}. Let
$\alpha$ be the solution of the elliptic sine-Gordon equation (ESGE) given by $\tan(\alpha/4)={\rm e}^{x_1}$. By solving the B\"acklund transformation~\eqref{BTcase6} for the parameters $\phi_1=0$, and any $\phi_2\in (0,\infty)$, we obtain $\alpha_1$ and $\alpha_2$ solutions of the ESHGE, given as follows:
\begin{figure}[t]
 \centering
 \begin{minipage}{.32\linewidth}
 \centering
 {\includegraphics[height=100pt]{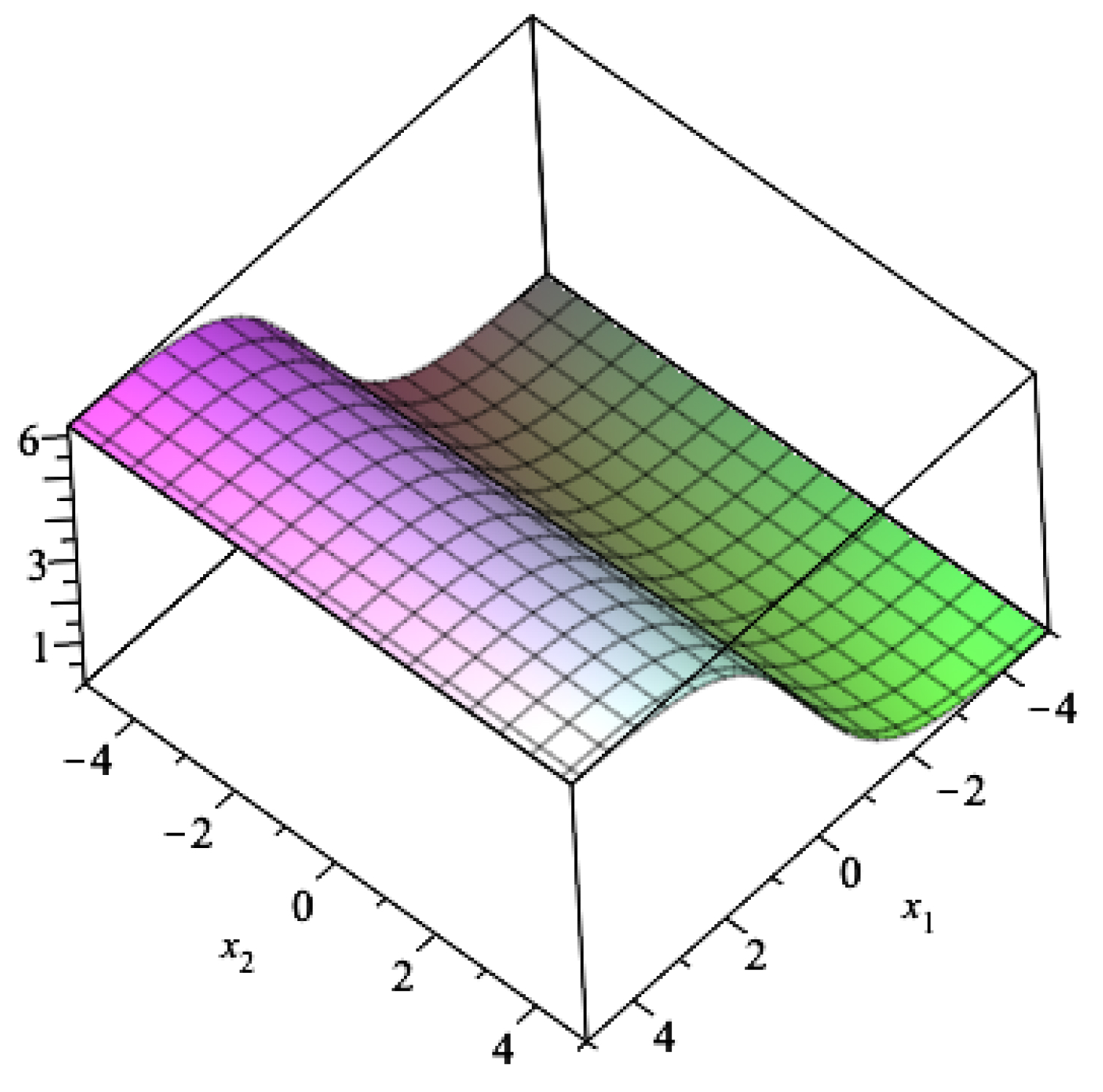}}\\
 $\alpha$
 \end{minipage}
 \begin{minipage}{.32\linewidth}
 \centering
 {\includegraphics[height=100pt]{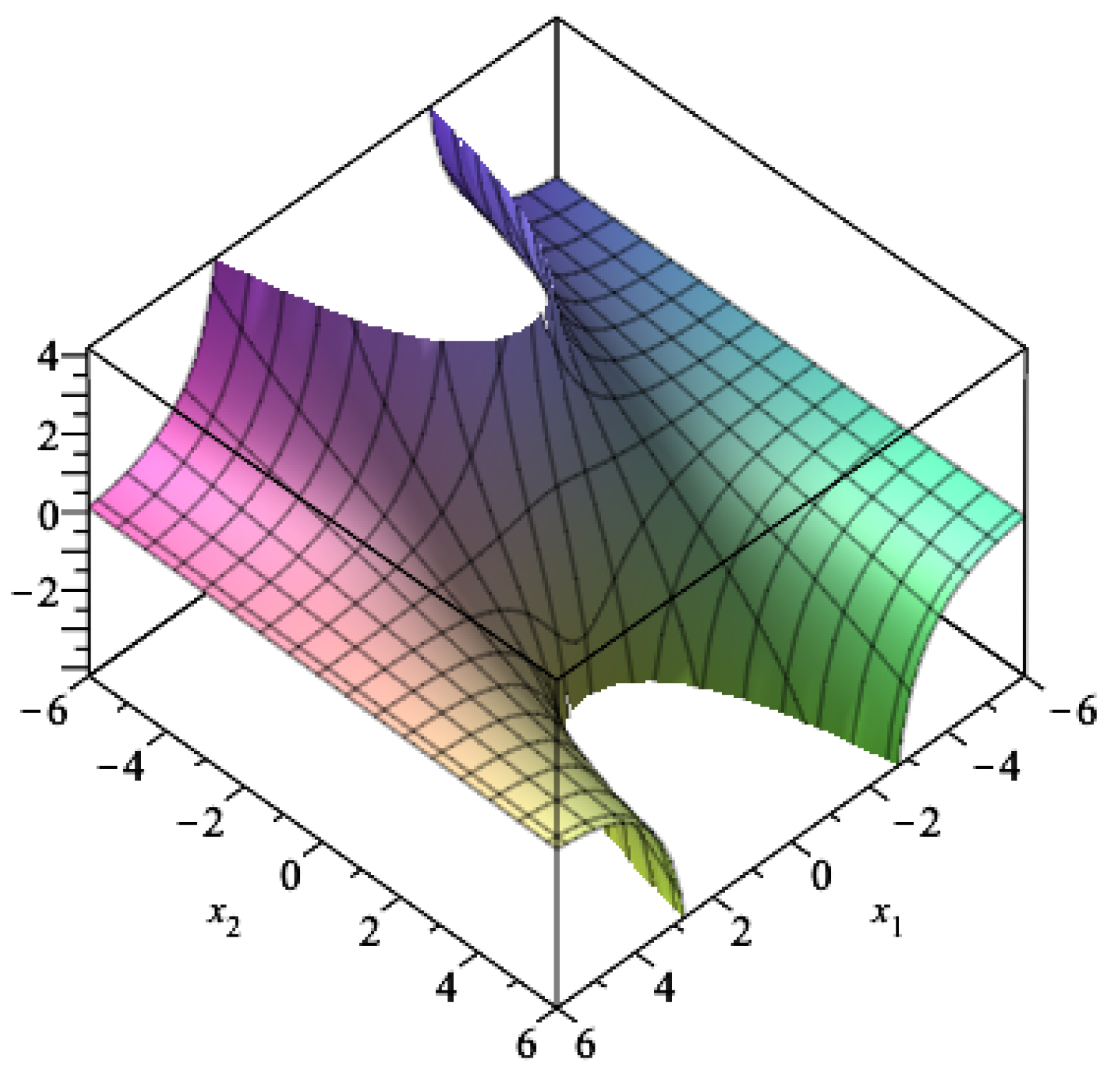}}\\
 $\alpha_1$ \\
\bigskip\bigskip\bigskip\bigskip
 {\includegraphics[height=100pt]{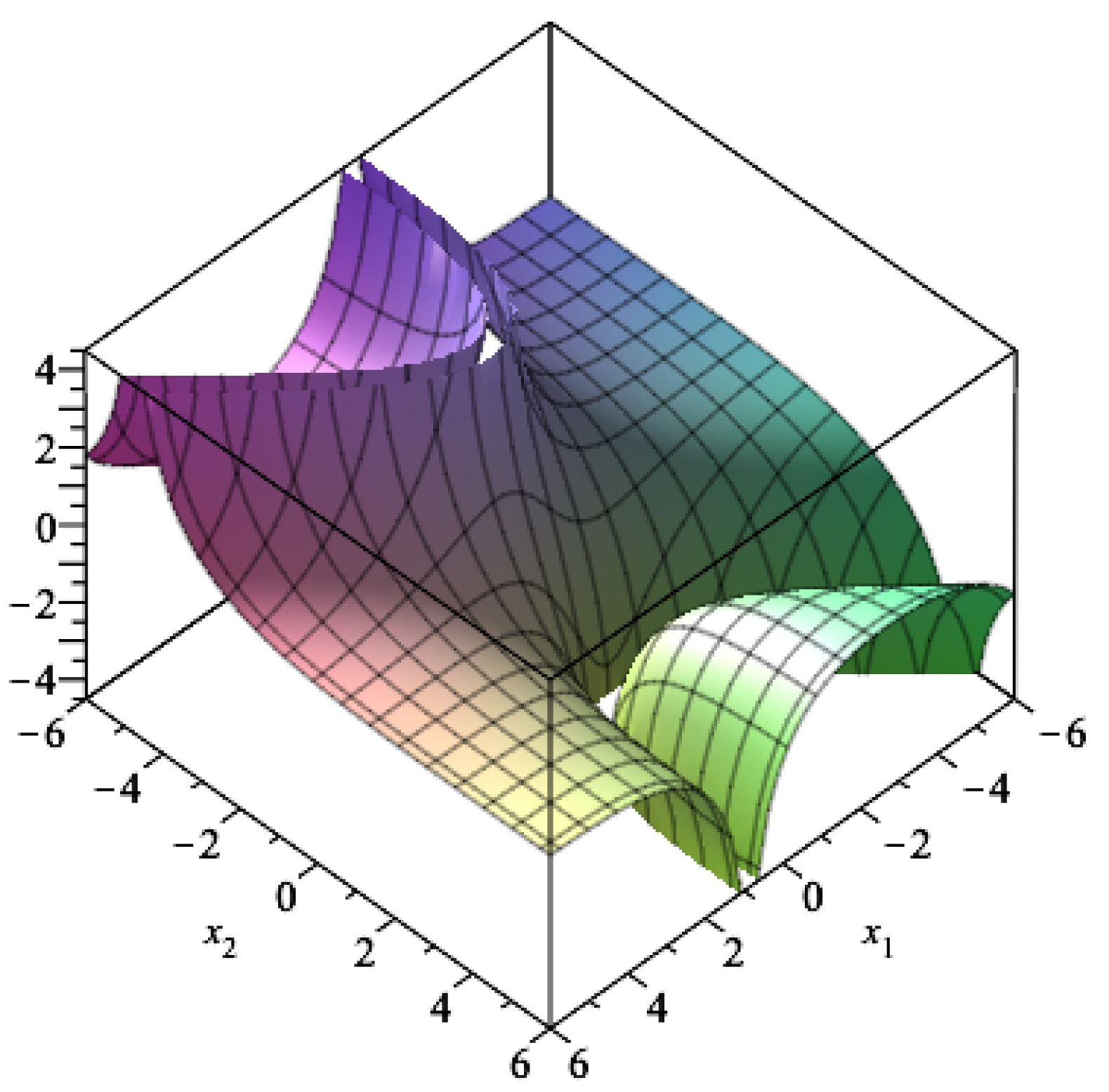}}\\
 $\alpha_2$
 \end{minipage}
 \begin{minipage}{.32\linewidth}
 \centering
 {\includegraphics[height=100pt]{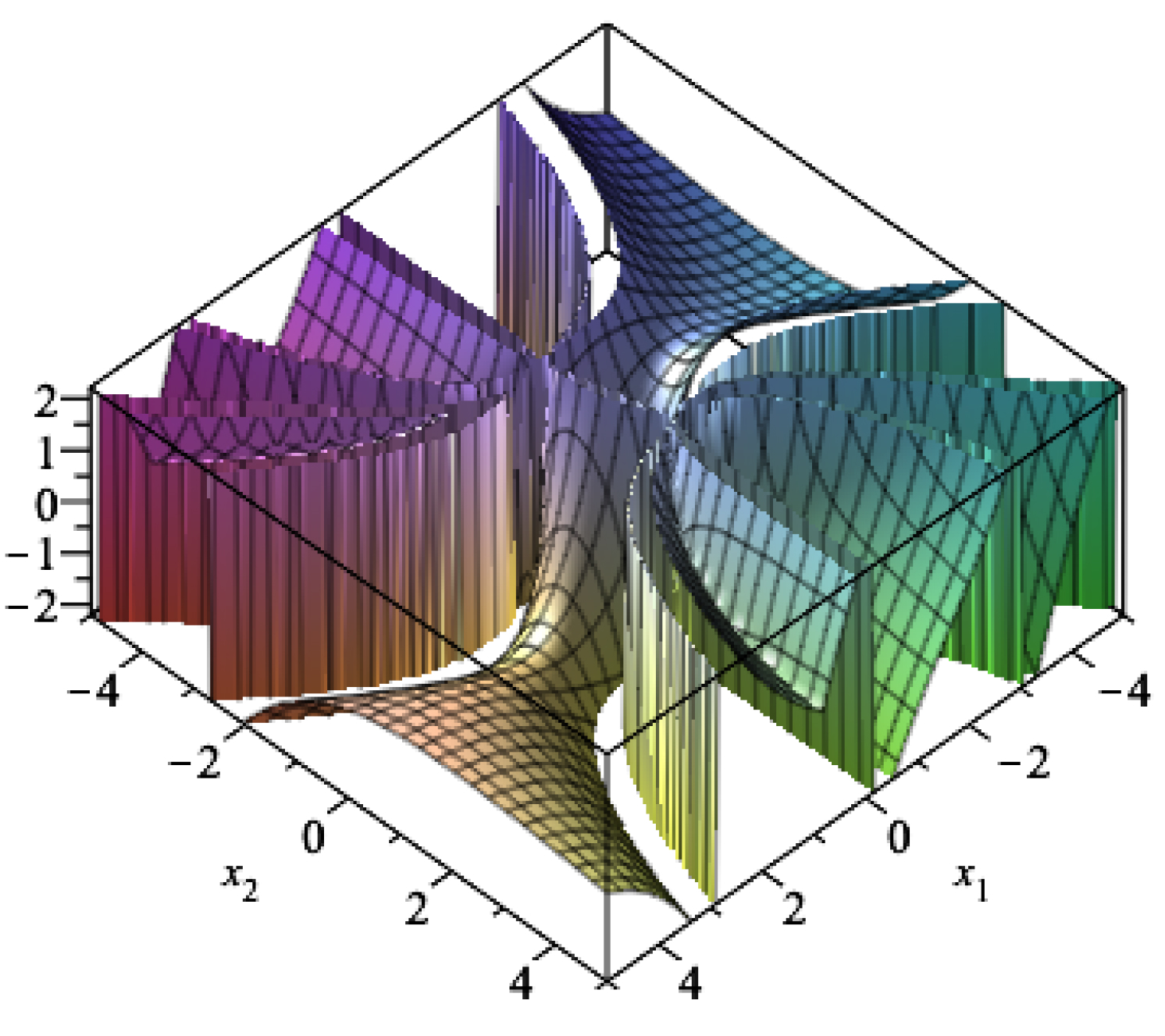}}\\
 $\alpha_{12}$
 \end{minipage}
 \caption{$\alpha$ is a solution of ESGE. $\alpha_1$ and $\alpha_2$ are solutions of the ESHGE obtained from $\alpha$ by the B\"acklund transformations with parameter $\phi_1=0$ and $\phi_2=\ln(\sqrt{2}+1)$ respectively. $\alpha_{12}$ is a solution of the ESGE obtained from $\alpha_1$ by using the superposition formula with parameters $\phi_1$ and $\phi_2$.}\label{fig:caso6alpha}
\end{figure}
\begin{gather}
\tanh\Bigl(\frac{\alpha_1}{4}\Bigr)=(-x_2+c_1)\sech x_1\label{alpha1caso6},\\
\tanh\Bigl(\frac{\alpha_2}{2}\Bigr)=\displaystyle{\frac{ \sinh{\phi_2 (\sinh \xi -\sinh x_1) }}
 {\sinh \xi \sinh x_1+1-\cosh\phi_2\cosh \xi \cosh x_1} } \label{alpha2caso6},
\end{gather}
where $\xi= (x_1+\sinh(\phi_2) x_2)/\cosh \phi_2 +c_2$ and $c_2\in{\mathbb{R}}$.
 The solution
$\alpha_1$ tends to $\pm\infty$ on two curves given by $x_2=\pm \cosh x_1 +c_1 $ (see Figure \ref{fig:caso6alpha}). The solution $\alpha_2$ is defined on the complement of two curves given by $\sinh \phi_2 (\sinh \xi -\sinh x_1) \pm (\sinh \xi \sinh x_1+1-\cosh\phi_2\cosh \xi \cosh x_1)=0$, where it tends to $\pm\infty$.
 The superposition formula \eqref{superposcase6}, defines explicitly a function $\alpha_{12}$ which is a new solution of the ESGE, given by
\[
\tan \Bigl(\frac{\alpha_{12}}{2}\Bigr)= \frac{\sinh x_1 P_{12}+Q_{12}}{\sinh x_1Q_{12}-P_{12}},
\]
where $P_{12}=\sinh\phi_2 \sinh((\alpha_2-\alpha_1)/2)$ and $Q_{12}=\cosh((\alpha_2-\alpha_1)/2)-\cosh\phi_2$ are determined by~\eqref{alpha1caso6} and~\eqref{alpha2caso6}. Observe that the solution
$\alpha_{12}$ tends to $\pm\infty$ on many curves defined by
\[
\arctan( (\sinh x_1 P_{12}+Q_{12})/(\sinh x_1Q_{12}-P_{12}))=\pi/2+k\pi,
\]
where $k$ is any integer (see Figure \ref{fig:caso6alpha}).
In particular, by choosing $\phi_2=\ln(1+\sqrt{2})$ and the constants $c_1=c_2=0$, one can visualize
the graph of the functions $\alpha$, $\alpha_1$, $\alpha_2$, $\alpha_{12}$ in Figure \ref{fig:caso6alpha}.
\end{Example}

Now we consider the surfaces in ${\mathbb{R}}^3_1$ of constant Gaussian curvature $K=-1$, that correspond to the functions above. The time-like surface given by
\begin{equation}\label{superfESGE}
X(x_1,x_2)=(x_1-\tanh x_1, \sech x_1\cosh x_2, \sech x_1\sinh x_2),
\end{equation}
corresponds to the solution $\alpha$. The geometric B\"acklund transformation of the surface $X$ for the parameters $\phi_1$ and $\phi_2$, provides two space-like surfaces of constant negative curvature
given by
\begin{gather*}
X_j=X+\cosh \phi_j(-\cosh(\alpha_j/2)/\tanh(x_1)X_{x_1}+ \sinh(\alpha_j/2)\cosh(x_1)X_{x_2}) ,\qquad j=1,2.
\end{gather*}
The fourth surface $X_{12}$ given by
\[
X_{12}=X_1+ \cosh\phi_2 \Bigl( 1-\tanh^2\Bigl(\frac{\alpha_1}{4}\Bigr) \Bigr) \biggl[ \dfrac{ \cos(\alpha_{12}/2) } {1+\tanh^2({\alpha_1}/{4})} X_{1,x_1}+ \dfrac{ \sin(\alpha_{12}/2) } {2 \tanh({\alpha_1}/{4})} X_{1,x_2} \biggr],
\]
where $ \tanh(\alpha_1/{4})$ is given by \eqref{alpha1caso6}, is a time-like surface that corresponds to $\alpha_{12}$ and it is associated to $X_1$ and to $X_2$ by a B\"acklund transformation with $\phi_2$ and $\phi_1$, respectively.

\subsection*{Acknowledgements}
Filipe Kelmer was partially supported by CNPq grant 132908/2015-8 and
CAPES/Brazil-Finance Code 001. Keti Tenenblat was partially supported by CNPq, Brasil, Proc. 311916/2021-0 and CAPES/Brazil-Finance Code 001.
 The authors would like to thank the anonymous referees whose comments gave a relevant contribution to the improvement of this paper.

\pdfbookmark[1]{References}{ref}
\LastPageEnding

\end{document}